\documentclass{amsart}
\usepackage[english]{babel}
\usepackage{latexsym}
\usepackage{amssymb}
\usepackage{amsmath}
\usepackage{amscd}
\usepackage[matrix,arrow,curve]{xy}
\usepackage{graphicx}
\usepackage[cp1251]{inputenc}
\CompileMatrices 
 \DeclareMathOperator{\TrueIm}{Im}
 \DeclareMathOperator{\Aut}{Aut}
 
 \DeclareMathOperator{\idd}{id}
 
\DeclareMathOperator{\repq}{Rep} \DeclareMathOperator{\HOMM}{Hom}
 \DeclareMathOperator{\Spec}{Spec}
 
\DeclareMathOperator{\lspan}{span} \DeclareMathOperator{\Grass}{Gr}
\DeclareMathOperator{\IHom}{IHom} 
\newtheorem{theor}{$\phantom{sss}$Theorem}
\newtheorem{utver}{$\phantom{sss}$Proposition}
\newtheorem{laemma}{$\phantom{sss}$Lemma}
\newtheorem{corol}{$\phantom{sss}$Corollary}

\begin{document}
\fontfamily{ptm} \fontsize{11pt}{14pt} \selectfont
\author[S. Fedotov]{Stanislav Fedotov}
\title[Framed moduli and Grassmannians]{Framed moduli and Grassmannians of submodules}
\thanks{Supported by grant RFFI 09-01-90416 - Ukr-f-a}
\address{Moscow State University, department of Higher Algebra
\endgraf Email: glwrath@yandex.ru}
\sloppy

\maketitle

\begin{abstract} In this work we study a realization of moduli
spaces of framed quiver representations as Grassmannians of
submodules devised by Marcus Reineke. Obtained is a generalization
of this construction for finite dimensional associative algebras and
for quivers with oriented cycles. As an application we get an
explicit realization of fibers for the moduli space bundle over the
categorical quotient for the quiver $A_{n-1}^{(1)}$.
\end{abstract}

\section{Introduction}

A quiver $Q$ is a diagram of arrows, determined by two finite sets
$Q_0$ (the set of ``vertices'') and $Q_1$ (the set of ``arrows'')
with two maps $h,t: Q_1\rightarrow Q_0$ which indicate the vertices
at the head and tail of each arrow. A representation $(W, \varphi)$
of $Q$ consists of a collection of finite dimensional $\Bbbk$-vector
spaces $W_i$, for each $i\in Q_0$, together with linear maps
$\varphi_a: W_{ta}\rightarrow W_{ha}$, for each $a\in Q_1$. The
dimension vector $\alpha\in\mathbb{Z}^{Q_0}$ of such a
representation is given by $\alpha_v = \dim_{\Bbbk}{W_i}$. A
morphism $f:(W_i, \varphi_a)\rightarrow (U_i,\psi_a)$ of
representations consists of linear maps $f_i: W_v\rightarrow U_v$,
for each $i\in Q_0$, such that $f_{ha}\varphi_a = \psi_af_{ta}$, for
each $a\in Q_1$. Evidently, it is an isomorphism if and only if each
$f_i$ is. Thus, isomorphism classes of representations of $Q$ with
dimension vector $\alpha$ coincide with orbits of the action of
$GL(\alpha) = \prod_{i\in Q_0}GL_i(\Bbbk)$ on the representation
space $\repq(Q, \alpha)$.

In studying quiver representations standard approaches of Invariant
Theory often fail because the algebra of invariants is poor or even
trivial, as in case of the quiver without oriented cycles, and so
the categorical quotient $\repq(Q, \alpha)/\!\!/GL(\alpha) :=
\Spec\Bbbk[\repq(Q, \alpha)]^{GL(\alpha)}$ is a point. Constructions
of Geometric Invariant Theory may help to compensate this defect.
Indeed, one can consider the trivial linearization twisted by a
character $\chi$ of $GL(\alpha)$, which restricts our attention to
an open subset of $\repq(Q, {\alpha})$, consisting of
$\chi$-semistable representations. Within the open set there are
more closed orbits and the corresponding algebraic quotient is more
interesting. In his paper~\cite{Kinge} A.D. King showed that the
notions of semistability and stability, that arise from Geometric
Invariant Theory, coincide with more algebraic notions, expressed in
the language of abelian categories. Namely, he devised a link
between this concept and the so-called $\theta$-stability. All the
characters of $GL({\alpha})$ are given by $\chi_{\theta}(g) =
\prod_{v\in Q_{0}}\det{(g_v)}^{\theta_v}$ for
$\theta\in\mathbb{Z}^{Q_0}$, and a representation $W$ is a
$\chi_{\theta}$-\linebreak(semi-)stable point of the variety
$\repq(Q,{\alpha})$ if and only if it is $\theta$-(semi-)stable as
an element of the abelian category $\repq(Q)$ (see~\cite[Section
2]{Kinge}). King proved the existence of coarse moduli space for
semistable representations of quivers and guaranteed that for stable
representations there is a fine moduli space. This technique allowed
a straightforward and convenient adaptation to the study of
representations of finite dimensional algebras~\cite[Section
4]{Kinge}.

An alternative approach to the problem considered was introduced by
B.~Huisgen-Zimmermann~\cite{BH-cl}. Let $A$ be a finite dimensional
associative algebra with unity. Fix a semisimple
$A$-module $T$ with projective cover $P$ and a positive
integer $m$. Denote by $\mathfrak{Grass}^T_m$ the
Grassmannian of all $(\dim P - m)$-dimensional submodules of the
radical $\mathrm{rad}{P}$. There is then a bijection between
$\Aut_AP$-orbits in $\mathfrak{Grass}^T_m$ and isomorphism classes
of $m$-dimensional modules with top isomorphic to $T$, sending $\Aut_AP\cdot C$ to
the isomorphism class of $P/C$. In~\cite{BH-cl} investigated are
such triples ($A$,~$T$,~$m$), that $\mathfrak{Grass}^T_m$ itself
provides a moduli space classifying the $d$-dimensional $A$-modules
with top $T$, up to isomorphism. It is also proved there that
$\mathfrak{Grass}^T_m$ admits an open covering by
representation-theoretically defined affine charts. For further
generalization and systematical treatment of the whole hierarchy of
moduli-parametrizing Grassmannian varieties see~\cite{BH-hier}.

Another possible way of applying Invariant Theory to the study of
quiver representations is to consider framed representations thus
achieving better precision at the expense of extending the
representation space. They first appeared in~\cite{Nakajima} as one
of the steps in the construction of Nakajima Varieties.

Let $Q$ be a quiver and $\alpha$ be a dimension vector. Fix an
additional dimension vector $\zeta$ and consider the space $\repq(Q,
\alpha, \zeta) := \repq(Q, \alpha)\oplus\bigoplus_{i\in
Q_0}\HOMM_{\Bbbk}(\Bbbk^{\alpha_i},\Bbbk^{\zeta_i})$. Its elements
are said to be {\it framed representations} of $Q$. Define a
$GL(\alpha)$-action on $\repq(Q, \alpha, \zeta)$ by $g\cdot(M,
(f_i)_{i = 1}^n) = (g\cdot M, (f_ig_i^{-1})_{i = 1}^n)$. A framed
representation $(M, f)$ is called {\it stable} if there is no
nonzero subrepresentation $N$ of $M$ which is contained in
$\ker{f}$. Denote by $\repq^s(Q, \alpha, \zeta)$ the space of stable
framed representations. One can show that the quotient of
$\repq^s(Q, \alpha, \zeta)$ is more efficient in orbits
discriminating than the standard categorical quotient. Furthermore,
it can be shown that it is in a sense reducible to King's
construction and thus enjoys all its properties.

Grassmanians of submodules of injective modules arise independently
in the course of this approach. In~\cite{Reineke} M. Reineke
obtained for acyclic quivers a realization of framed moduli space as
a Grassmannian of subrepresentations of an injective representation
depending only on dimension vectors $\alpha$ and $\zeta$. He further
investigated its cohomology (see also~\cite{Cald-Rein}) and
applications to quantum groups.

The aim of this work is to explore possible generalizations of
Reineke's construction. So far, there are two of them. First, we can
adapt it for quivers with relations, i.e. at least for finite
dimensional algebras. We may also try to eliminate the condition of
having no oriented cycles. Both possibilities are discussed below.

Section 2 is devoted to giving some basic results on the connection
between finite dimensional algebras and quivers. In Subsections 3.1
and 3.2 we remind the concept of stability as it was formulated by
A. D. King and A. N. Rudakov for abelian categories. We also introduce the
notion of framed representation space for finite dimensional
algebras.

A straightforward generalization of Reineke's ideas for finite
dimensional algebras is presented in Subsection 3.3. We prove that
for a finite dimension algebra and two dimensional vectors
$\alpha$ and $\zeta$ the quotient space $\repq^s(A, \alpha,
\zeta)/\!\!/GL(\alpha)$ is isomorphic to the
Grassmannian of submodules of a certain injective $A$-module $J$.

In Sections 4 to 6 we investigate a possible way to generalize this
result to quivers with oriented cycles. It is crucial for Reineke's construction
that $A$ is finite dimensional, since without it we are no
longer able to embed all stable framed representations in a finite
dimensional module. In general situation we thus come to considering
Grassmannians of submodules in infinite dimensional spaces, that
makes impossible to apply the usual techniques. So, having established a bijection between
such a Grassmannian and the moduli space we provide the former with
a structure of algebraic variety, but this gives us no information about the moduli
space itself. Hence we are forced to restrict
our attention to fibers of the moduli space bundle over the
categorical quotient.

From that point we work over the ground field $\mathbb{C}$. In sections 4 and 5 we study
the quiver $Q = A^{(1)}_{n-1}$ with cyclic orientation.
We obtain an embedding of
$\repq^s(Q, \alpha, \zeta)$ into an infinite dimensional
representation $J^{!!}$, with each $J^{!!}_i$ a space of holomorphic
vector functions on $\mathbb{A}^1$. Recall that by Procesi-Razmyslov's Theorem~\cite[Theorem 1]{LBP}
the algebra of invariant
polynomial functions on $\repq(Q, \alpha, \zeta)$ is generated by traces of
oriented cycles in $Q$. Since in $Q$ there is only one cycle $\tau_i$ of minimal length
starting at each vertex $i\in Q_0$, traces of $\tau_i^k$ are polynomials in
coefficients of characteristic polynomials $\chi_i$ of $\tau_i$.
Hence, the categorical
quotient $\repq(Q, \alpha, \zeta)/\!\!/GL(\alpha)$
 may be identified
with a subspace in the Cartesian product $\mathbb{C}[t]^n$, each point presented by a tuple
of characteristic polynomials. Let
$\overline{\chi} = (\chi_1,\ldots,\chi_n)$ be in $\repq(Q, \alpha,
\zeta)/\!\!/GL(\alpha)$, $\lambda_1,\ldots,\lambda_N$ be all
different roots of $\chi_1,\ldots,\chi_n$, and $r_{ij}$ be the
multiplicity of $\lambda_j$ as a root of $\chi_i$. For each $j =
1,\ldots,N$ consider a new dimension vector $\underline{r}_j$ with
$(\underline{r}_j)_i = r_{ij}$. Define submodules $J(\lambda_j,
\underline{r}_j)\subseteq J^{!!}$ by $J(\lambda_j,
\underline{r}_j)_i = \{\overline{g}\in J^{!!}\mid(\frac{d}{dt} -
\lambda_j\idd)^{\max_i{r_{ij}}}\cdot\overline{g} = 0\}$. These
representations are finite dimensional and it may be shown that the
fiber of the framed moduli space over
$(\chi_1,\ldots,\chi_n)\in\mathbb{C}[t]^n$ is isomorphic to the
product of the Grassmannians of their subrepresentations. The
details are discussed in Section 4. In Section 5 obtained is an
explicit presentation by equations in projective space of fibers of
moduli space bundle over the categorical quotient. Section 6 is
devoted to applying of this technique to any quiver $Q$ where all
oriented cycles pairwise commute. We consider the natural projection
$\pi_s:\repq^s(Q, \alpha, \zeta)/\!\!/GL(\alpha)\rightarrow\repq(Q,
\alpha, \zeta)/\!\!/GL(\alpha)$ and describe the fibers
$\pi_s^{-1}(x)$, for $x\in \repq(Q, \alpha, \zeta)/\!\!/GL(\alpha)$.
The main result here is that for every $x$ and two dimension vectors
$\alpha$ and $\zeta$ there is a quiver $Q^{\spadesuit}$, a dimension
vector $\widetilde{\alpha}\in(\mathbb{Z}_{\geqslant
0})^{Q^{\spadesuit}_0}$, and a finite dimensional representation
$W^{\spadesuit}$ of $Q^{\spadesuit}$ such that the fiber
$\pi_s^{-1}(x)$ is isomorphic to the Grassmannian of $\widetilde{\alpha}$-dimensional submodules of
$W^{\spadesuit}$.

After having written this paper the author found that fibers of the projection
$\pi_s$ were described by J. Engel and M. Reineke for arbitrary quivers using Luna's
stratification as nilpotent parts of the framed moduli space for some new
quiver $Q'$ and dimension vectors $\alpha'$ and $\zeta'$ (see~\cite[Theorem 4.1]{smooth}). Theorem 3 together
with Proposition 4 thus imply that each fiber of $\pi_s$ is isomorphic to the Grassmannian
of submodules of a module over some finite dimensional algebra with a certain
dimension vector. Our description is less general than the one given in~\cite{smooth}, but more explicit
and convenient when applicable. 

The author thanks Ivan V. Arzhantsev and Markus Reineke,
who also brought to his attention the paper~\cite{smooth} and the universal description
of layers of $\pi_s$, for useful discussions.

\section{Background information}

First, recall some general facts concerning finite dimensional
algebras and their connection with quivers.

Let $A$ be an associative finite dimensional algebra over an
arbitrary field $\Bbbk$. An element $e\in A$ is called {\it
idempotent} if $e^2 = e$. We say that two idempotents $e_1$ and
$e_2$ are {\it orthogonal} if $e_1e_2 = e_2e_1 = 0$. An idempotent
$e$ is {\it primitive} if it is not a sum of two nonzero orthogonal
idempotents. It is well known that for an algebra $A$ with unity
there always exists a decomposition $1 = e_1 + \ldots + e_n$, where
$e_i$ are primitive orthogonal idempotents. Note that this
decomposition induces a decomposition $A = Ae_1\oplus\ldots\oplus
Ae_n$ of the regular $A$-module called {\it Peirce decomposition}.

An algebra is said to be {\it splitting} if the quotient
$A/\mathfrak{r}$, where $\mathfrak{r}$ stands for the radical of
$A$, is isomorphic to a direct product of matrix algebras over the
ground field. Note that all algebras over an algebraically closed
field are splitting. An algebra is called {\it basic} if
$A/\mathfrak{r}$ is isomorphic to a direct product of division
rings. This condition is equivalent to the following: in the decomposition $A =
P_1\oplus\ldots\oplus P_k$ of the regular module, where all $P_i$
are indecomposable projective modules, all the summands are pairwise
non-isomorphic,~\cite[Theorem 3.5.4]{Drozd}.

Fix a decomposition $1 = e_1 + \ldots + e_n$ of the unity $1\in A$,
where all $e_i$ are primitive orthogonal idempotents. It is not hard
to see that every $A$-module $M$ as a $\Bbbk$-vector space may be
decomposed as $M = \bigoplus_{i = 1}^nM_i$, where $M_i = e_iM$. The
{\it dimension vector} of $M$ is the vector $\alpha =
\underline{\dim}{M}$ with $\alpha_i = \dim{M_i}$. Decomposing in
this way the ideals $Ae_i$ of $A$ (which are submodules of the
regular module), we obtain the {\it two-sided Peirce decomposition}
of $A$: $A = \bigoplus_{i,j}e_jAe_i$. The components $e_jAe_i$ are
neither left, nor right ideals, but they provide a convenient matrix
interpretation of the elements of $A$ (see~\cite[Chapter 1, \S
7]{Drozd}).

The set of all $A$-modules with dimension vector $\alpha$ will be
denoted by $\repq(A,\alpha)$. The group $GL(\alpha) = \prod_{i =
1}^nGL_{\alpha_i}(\Bbbk)$ acts naturally on this set, each factor
acting by base change in $M_i$. Namely, for an element $a\in
e_jAe_i$ the corresponding operator $\varphi(a)$ of the
representation $A\rightarrow L(M)$ maps $M_i$ to $M_j$ and all the
rest components to zero; thus we may define the action $g = (g_t)_{t
= 1}^n\in\prod_{i = 1}^nGL(M_i) = GL(\alpha)$ as follows:
$(g\cdot\varphi(a))(m) = (g_j\varphi(a)g_i^{-1})(m)$, for all $g\in
G$, $m\in M_i$. Since $A$ admits the two-sided Peirce decomposition
$A = \bigoplus_{i,j}e_jAe_i$, the actions is well defined.

Now we remind the connection between algebras and quivers. Let $\Bbbk$ be a field.
For a quiver $Q$ one defines a {\it path algebra} $\Bbbk Q$. As a
linear space it is the span of all paths in $Q$, including those of
length $0$, which we identify with vertices of $Q$. Multiplication
in $\Bbbk Q$ is defined by
$$\sigma\cdot\tau = \begin{cases} \sigma\tau,\,\mbox{if it is a path in
$Q$,}\\ 0,\,\mbox{otherwise,}
\end{cases}$$
for two paths $\sigma$ and $\tau$ in $Q$.

A {\it relation} in $Q$ is a $\Bbbk$-linear combination
of paths in $Q$ of length not less than 2 with the same source and
target. For a set of relations $\rho$ denote by $\langle\rho\rangle$
the ideal of the algebra $\Bbbk Q$ generated by these relations.

\begin{theor}~\cite[Prop. II.2.5]{A-R} For a finite dimensional algebra $A$
the category of finitely generated $A$-modules is equivalent to the
category of finitely generated $\Gamma$-modules for some basic
algebra $\Gamma$.
\end{theor}

Thus the problem of classifying the representations of arbitrary
finite dimensional algebras can be in a sense reduced to the case of
basic algebras. So we will be considering only basic algebras.

We now remind briefly the procedure of assigning a quiver $Q(A)$ to
a splitting basic algebra $A$. Let $\mathfrak{r}$ be the radical of
$A$, $1 = \overline{e}_1 + \ldots + \overline{e}_n$ be a
decomposition of unity in $A/\mathfrak{r}$, $1 = e_1 + \ldots + e_n$
the corresponding decomposition of unity in $A$; further, let $W =
(\mathfrak{r}/\mathfrak{r}^2)$. Denote $t_{ij} =
\dim_{\Bbbk}{e_jWe_i}$. Now set $Q(A) = (Q_0, Q_1)$ with $Q_0 =
\{1,\ldots,n\}$ and $t_{ij}$ arrows from the $i$-th vertex to the
$j$-th one.

Let $(\Bbbk Q)_{\geqslant1}$ be the ideal generated by the arrows in
$Q$. An ideal $I\vartriangleleft\Bbbk Q$ is said to be {\it regular}
if $(\Bbbk Q)_{\geqslant1}^2\supseteq I\supseteq(\Bbbk
Q)_{\geqslant1}^t$, for some $t\geqslant2$.

\begin{theor}~\cite[Theorem III.1.9]{A-R} Every splitting basic finite dimensional algebra
with quiver $Q$ is isomorphic to a factor algebra $\Bbbk Q / I$,
where $I$ is a regular ideal.
\end{theor}

\begin{corol} For a splitting basic finite dimensional algebra
$A$ there is a set of relations $\rho$ such that $A\cong\Bbbk
Q(A)/\langle\rho\rangle$.
\end{corol}

On the level of representation spaces this correspondence looks as
follows: the set $\repq(A,\alpha)$ is a Zariski closed subvariety of
$\repq(Q(A),\alpha)$, since it is the subset where $X_p \equiv 0$
for all $p\in \langle\rho\rangle$ (for a representation
$X\in\repq(Q(A),\alpha)$ and an element $\lambda_1a_{i_{11}}\ldots
a_{i_{1k(1)}} + \ldots + \lambda_sa_{i_{s1}}\ldots a_{i_{sk(s)}}$ we
denote by $X_p$ the linear transformation
$\lambda_1X_{a_{i_{11}}}\ldots X_{a_{i_{1k(1)}}} + \ldots +
\lambda_sX_{a_{i_{s1}}}\ldots X_{a_{i_{sk(s)}}}$). Sometimes we will
denote this subvariety by $\repq(Q(A), \rho, \alpha)$.

\section{Framed representations of finite dimensional algebras}

\subsection{Semistabile representations}
Consider a quiver $Q$. A {\it character} of the category $\repq(Q)$
of representations of $Q$ is a linear function $\theta: \mathbb{Z}
Q_0\rightarrow\mathbb{Z}$ (in other words, to each vertex of the
quiver this function assigns an integer). For $X\in\repq(Q)$ define
$\theta(X) = \theta(\underline{\dim}{X})$. A representation $X$ is
{\it $\theta$-semistable} (respectively {\it $\theta$-stable}) if
$\theta(X) = 0$ and $\theta(Y)\geqslant 0$ for all
subrepresentations $Y\subseteq X$ (respectively $\theta(Y)
> 0$). This approach devised by A. D. King, was generalized and reformulated in a
more flexible form by A. N. Rudakov~\cite[\S 3]{Rudak}.

Consider two characters $\theta, \kappa: \mathbb{Z}
Q_0\rightarrow\mathbb{Z}$, such that $\kappa(d)\geqslant 0$, for
every vector $d$ with nonnegative components (i.e. for every
dimension vector). Define a {\it slope}
$\mu:\repq(Q)\backslash\{0\}\rightarrow\mathbb{Q}$ by $\mu(X) =
\mu(\underline{\dim}{X}) =
\frac{\theta(\underline{\dim}{X})}{\kappa(\underline{\dim}{X})}$. A
representation $X$ is called {\it $\mu$-semistable} (respectively
{\it $\mu$-stable}) if $\mu(Y)\leqslant \mu(X)$ (respectively
$\mu(Y) < \mu(X)$) for each proper nonzero subrepresentation $Y$ of
$X$. Denote by $\repq^{ss}_{\mu}(Q, \alpha)$ (respectively
$\repq^{s}_{\mu}(Q, \alpha)$) the set of $\mu$-semistable
(respectively $\mu$-stable) representations of $Q$ with dimension
vector $\alpha$. The following simple lemma shows the connection
between this notion and King's construction.

\begin{laemma} Let $\alpha$ be a dimension vector and $\mu$ be a slope.
There is a character $\xi$ (depending on $\alpha$) such that
$\repq^{ss}_{\mu}(Q, \alpha) = \repq^{ss}_{\xi}(Q, \alpha)$ and
$\repq^{s}_{\mu}(Q, \alpha) = \repq^{s}_{\xi}(Q, \alpha)$.
\end{laemma}

\begin{proof} Consider $\mu = \frac{\theta}{\kappa}$, where $\kappa$
is a linear function taking nonnegative values on dimension vectors.
Set $\xi(d) = \mu(\alpha)\kappa(d) - \theta(d)$. We know that
$X\in\repq^{ss}_{\mu}(Q, \alpha)$ if and only if
$\mu(Y)\leqslant\mu(X)$ for every non-trivial subrepresentation $Y$,
that is if and only if
$\frac{\theta(Y)}{\kappa(Y)}\leqslant\mu(\alpha)$. This means that
$0\leqslant\mu(\alpha)\kappa(Y) - \theta(Y) = \xi(Y)$. Since $\xi(X)
= \xi(\alpha) = 0$, we obtain that $X\in\repq^{ss}_{\mu}(Q, \alpha)
\Leftrightarrow X\in\repq^{ss}_{\xi}(Q, \alpha)$. The second part of
the lemma is proved using the same arguments.
\end{proof}

From now on, since for the subsets $\repq^{ss}_{\xi}(Q,\alpha)$ and
$\repq^{s}_{\xi}(Q,\alpha)$ the existence of the categorical
quotient was proved in~\cite{Kinge}, we can use the notion of moduli
spaces of $\mu$-(semi-)stable points, where $\mu$ is a slope.
Namely, denote by $\mathcal{M}^{ss}_{\mu}(Q, \alpha)$ (respectively
by $\mathcal{M}^{s}_{\mu}(Q, \alpha)$) the categorical quotient
$\repq^{ss}_{\mu}(Q, \alpha)/\!\!/GL(\alpha)$ (respectively
$\repq^{s}_{\mu}(Q, \alpha)/\!\!/GL(\alpha)$)

Now let $A$ be a finite dimensional algebra, $Q = Q(A)$ be its
quiver, $\rho$ be a set of relations such that $A\cong\Bbbk
Q/\langle\rho\rangle$, and $\mu: \mathbb{Z}Q_0\rightarrow\mathbb{Q}$
be a slope. The fact that the set $\repq(A,\alpha)$ is embedded in
$\repq(Q, \alpha)$ as a Zariski closed $GL(\alpha)$-invariant
subvariety allows us to define $\mu$-semistable and $\mu$-stable
$A$-modules and, consequently, the subsets $\repq^{ss}_{\mu}(A,
\alpha)$ and $\repq^{ss}_{\mu}(A, \alpha)$ and their categorical
quotients $\mathcal{M}^{ss}_{\mu}(A, \alpha) := \repq^{ss}_{\mu}(A,
\alpha)/\!\!/GL(\alpha)$ and $\mathcal{M}^{s}_{\mu}(A, \alpha) :=
\repq^{s}_{\mu}(A, \alpha)/GL(\alpha)$ (the latter turns out to be a
geometric quotient, so there is a single fraction bar).

\subsection{Framed representations}

Let $A$ be a (splitting basic) finite dimensional algebra with
$|Q(A)_0| = n$. Fix two dimension vectors $\alpha,
\zeta\in(\mathbb{Z}_{\geqslant0})^n$ and consider the extended
representation space $\repq(A,\alpha,\zeta) = \{(M,f)\mid
M\in\repq(A,\alpha), f = (f_i)_{i = 1}^n: M_1\oplus\ldots\oplus
M_n\rightarrow V_1\oplus\ldots\oplus V_n \mbox{ is a graded linear
map}\} \cong
\repq(A,\alpha)\oplus\bigoplus_{i=1}^n\HOMM_{\Bbbk}(\Bbbk^{\alpha_i},
V_i)$, where $\dim{V_i} = \zeta_i$. Define the action of
$GL(\alpha)$ on this space by $g\cdot(M, (f_i)_{i = 1}^n) = (g\cdot
M, (f_ig_i^{-1})_{i = 1}^n)$. Elements of $\repq(A,\alpha,\zeta)$
will be called {\it framed representations} of the algebra $A$.

{\bf Definition.} A pair $(M, f)\in\repq(A,\alpha,\zeta)$ is {\it
stable} if no nonzero submodule of $M$ is contained in $\ker{f}$.
The set consisting of such pairs will be denoted by
$\repq^s(A,\alpha,\zeta)$.

Let $\zeta$ be a dimension vector. We introduce a new quiver
$\widetilde{Q}$ with $\widetilde{Q}_0 = Q_0\cup\{\infty\}$, the
arrow of $\widetilde{Q}$ being those of $Q$ together with $\zeta_i$
arrows from $i$ ($i\in Q_0$) to $\infty$. We also extend the
dimension vector $\alpha$ to $\widetilde{\alpha}$, setting
$\widetilde{\alpha}_i = \alpha_i$ for $i = 1,\ldots,n$ and
$\widetilde{\alpha}_{\infty} = 1$.

Observe that the elements of $\rho$ are relations in
$\widetilde{Q}$; consider the ideal $I = \langle\rho\rangle$ in
$\Bbbk\widetilde{Q}$. Then $(\Bbbk
\widetilde{Q})_{\geqslant1}^2\supseteq I\supseteq(\Bbbk
\widetilde{Q})_{\geqslant1}^{t + 1}$ for $t\geqslant2$ such that $
\langle\rho\rangle_{\Bbbk Q}\supseteq(\Bbbk Q)_{\geqslant1}^t$ (here
$\langle\rho\rangle_{\Bbbk Q}$ stands for the ideal of $\Bbbk Q$
generated by $\rho$). The last statement is not that obvious; but
recall that all the new arrows terminate in $\infty$, which means
that no path starts in this vertex. Therefore, if $p$ is a path of
length not less than $t + 1$ in $\widetilde{Q}$, then either it is
entirely contained in $Q$ and so $p\in I$, since $
\langle\rho\rangle\supseteq(\Bbbk Q)_{\geqslant1}^t = (\Bbbk
Q)_{\geqslant t}$, or $p = bq$, where $q$ is a path entirely
contained in $Q$ and $b$ is an arrow ending in $\infty$. But in the
latter case the length of $q$ is not less than $t$ yielding that
$q\in I$. Thus, the algebra $\widetilde{A} = \Bbbk\widetilde{Q}/I$
is finite dimensional.

Further, for $\beta\in\mathbb{Z}\widetilde{Q}_0$ set $\theta(\beta)
= - \beta_{\infty}$, $\kappa(\beta) = \sum_i\beta_i$ and consider
the corresponding slope $\mu = \frac{\theta}{\kappa}$. We are now
going to use Corollary 1 to get an interpretation of our notion of
stability in the spirit of King's construction.

\begin{utver} The sets
$\repq^s(A,\alpha,\zeta)$ and $\repq^{s}_{\mu}(\widetilde{A},
\widetilde{\alpha})$ are isomorphic as algebraic varieties.
\end{utver}

\begin{proof} The isomorphism at the level of quivers, i. e. for $A =
\Bbbk Q$, is proved in~\cite[Prop. 3.3]{Reineke}. We just point out
that crucial here is the existence of the $GL(\alpha)$-invariant
isomorphism $\HOMM_{\Bbbk}(\Bbbk^{\alpha_i},\Bbbk^{\zeta_i})\cong
\HOMM_{\Bbbk}(\Bbbk^{\alpha_i},\Bbbk)^{\zeta_i}$ providing
$$\repq(Q,
\alpha,\zeta)\cong\bigoplus_{a:i\rightarrow
j}\HOMM(\Bbbk^{\alpha_i},\Bbbk^{\alpha_j})\oplus\bigoplus_{i\in
Q_0}\HOMM_{\Bbbk}(\Bbbk^{\alpha_i},\Bbbk^{\zeta_i})\cong$$
$$\cong\bigoplus_{a:i\rightarrow
j}\HOMM(\Bbbk^{\alpha_i},\Bbbk^{\alpha_j})\oplus\bigoplus_{i\in
Q_0}\HOMM_{\Bbbk}(\Bbbk^{\alpha_i},\Bbbk)^{\zeta_i}\cong\repq(\widetilde{Q},\widetilde{\alpha}).
\eqno(1)$$

In order to pass to the general case we need to show that the image
of $\repq(A,\alpha, \zeta) = \left\{(X, f)\in\repq(Q, \alpha,
\zeta)\mid X_p = 0\,\forall p\in\rho\right\}$ under (1) is
$\repq(\widetilde{A}, \widetilde{\alpha}) =
\left\{Y\in\repq(\widetilde{Q}, \widetilde{\alpha})\mid Y_p =
0\,\forall p\in\rho\right\}$. But this follows from the fact that
the relations $\rho$ only affect the summand
$\bigoplus_{a:i\rightarrow
j}\HOMM(\Bbbk^{\alpha_i},\Bbbk^{\alpha_j})$, which is common for
both sets $\repq(Q, \alpha, \zeta)$ and $\repq(\widetilde{Q},
\widetilde{\alpha})$, while the restriction of (1) to this summand
is the identity map.

Now it is left to note that $\repq^s(A, \alpha,\zeta) = \repq^s(Q,
\alpha, \zeta)\cap\repq(A,\alpha,\zeta)$ and
$\repq^s_{\mu}(\widetilde{A}, \widetilde{\alpha}) =
\repq^s_{\mu}(\widetilde{Q},
\widetilde{\alpha})\cap\repq(\widetilde{A},\widetilde{\alpha})$,
finishing the proof.
\end{proof}

\begin{corol} The moduli spaces $\mathcal{M}^s(A,\alpha,\zeta):= \repq^s(A,\alpha,\zeta)/\!\!/
GL(\alpha)$
 and $\mathcal{M}^s_{\mu}(\widetilde{A},\widetilde{\alpha})$ are isomorphic.
\end{corol}

\begin{proof} It is a straightforward consequence of $GL(\alpha)$-invariance of the isomorphism (1).
\end{proof}

\begin{corol} The quotient $\repq^s(A,\alpha,\zeta)\rightarrow\mathcal{M}^s(A,\alpha,\zeta)$ is geometric.
Moreover, if $\repq^s(A,\alpha,\zeta)$ is nonempty, then
$\mathcal{M}^s(A,\alpha,\zeta)$ is a smooth projective variety.
\end{corol}

\begin{proof} We will prove these properties for
$\mathcal{M}^s_{\mu}(\widetilde{A},\widetilde{\alpha})$. First of
all, note that orbits of points from $\repq^s_{\mu}(\widetilde{A},
\widetilde{\alpha})$ are closed in $\repq^{ss}_{\mu}(\widetilde{A},
\widetilde{\alpha}) = \repq^{s}_{\mu}(\widetilde{A},
\widetilde{\alpha})$ (see the geometric definition of stability
in~\cite{Kinge}) and so the quotient is geometric. Furthermore,
standard results of algebraic geometry imply that
$\mathcal{M}^s_{\mu}(\widetilde{A},\widetilde{\alpha}) =
\mathcal{M}^{ss}_{\mu}(\widetilde{A},\widetilde{\alpha})$ is
projective over the categorical quotient $\mathcal{M}(\widetilde{A},
\widetilde{\alpha}) := \repq(\widetilde{A},
\widetilde{\alpha})/\!\!/GL(\alpha)$. Consider
$\Bbbk[\repq(\widetilde{A}, \widetilde{\alpha})]^{GL(\alpha)}$ as a
subalgebra in $\Bbbk[\repq(\widetilde{Q}, \widetilde{\alpha})]$. By
Procesi-Razmyslov's Theorem~\cite[Theorem 1]{LBP} the algebra
$\Bbbk[\repq(\widetilde{Q}, \widetilde{\alpha})]$ is generated by
traces of all oriented cycles in $\widetilde{Q}$. But since there is $t\geqslant 2$
such that $\langle\rho\rangle\supseteq
\Bbbk Q_{\geqslant t}$, some positive power of each oriented cycle
in $\widetilde{Q}$ lies in the ideal $I$ of $\Bbbk\widetilde{Q}$
generated by $\rho$, and hence in each representation belonging to
$\repq(\widetilde{A},\widetilde{\alpha}) =
\repq(\widetilde{Q},\rho,\widetilde{\alpha})$ all the oriented
cycles are nilpotent operators, and consequently their traces are
equal to zero. Thus, $\Bbbk[\repq(\widetilde{A},
\widetilde{\alpha})]^{GL(\alpha)} = \Bbbk$ implying that
$\mathcal{M}(\widetilde{A},
\widetilde{\alpha})\cong\left\{pt\right\}$ and, therefore,
$\mathcal{M}^s_{\mu}(\widetilde{A},\widetilde{\alpha})$ is a
projective variety.
\end{proof}

\subsection{The construction of the quotient space}

$\phantom{AAA}$

From now on for $\sigma\in\Bbbk Q$ we will denote by
$\overline{\sigma}$ its image $\sigma +
\langle\rho\rangle$ in $\Bbbk Q/\langle\rho\rangle\cong A$. It is easy
to see that in $\Bbbk Q/\langle\rho\rangle$ there is a (finite)
$\Bbbk$-basis $\Xi$ consisting of images of paths in $\Bbbk Q$. Its
elements will be referred to as paths in the algebra $A$. Denote by
$I_i$ the injective $A$-module associated with the $i$-th vertex of
the quiver. Recall that the corresponding representation from
$\repq(Q,\rho)$ may be described as follows: $(I_i)_j =
\lspan\left\{\overline{\tau}: j\rightsquigarrow i\right\}^*$, where
``$\overline{\tau}: j\rightsquigarrow i$'' means that
$\overline{\tau}$ is the image of a path $\tau$ starting in the
$j$-th vertex and ending in the $i$-th one; in this case $((I_i)_{a:
k\rightarrow l}f)(\overline{\tau}) =
f(\overline{\tau}\,\overline{a})$, where $\overline{\tau}:
l\rightsquigarrow i$. Consider the injective module $J :=
\bigoplus_{i\in Q_0}I_i\otimes_{\Bbbk} V_i$. Observe that as a
$\Bbbk$-linear space
$$J_i = e_iJ\cong\bigoplus_{j\in Q_0}(I_j)_i\otimes_{\Bbbk}V_j\cong
\bigoplus_{j\in
Q_0}\bigoplus_{\Xi\ni\overline{\tau}:i\rightsquigarrow
j}V_j\cong\bigoplus_{\Xi\ni\overline{\tau}:i\rightsquigarrow j}V_j$$

Now, given a point $(M, f)\in\repq(A,\alpha,\zeta)$ define a map
$\Phi_{(M,f)} = (\varphi_i)_{i\in Q_0}: M\rightarrow J$ by the
following rule:
$$\varphi_i = \bigoplus_{\Xi\ni\overline{\tau}: i\rightsquigarrow j}
f_j\overline{\tau}:
M_i\rightarrow\bigoplus_{\Xi\ni\overline{\tau}:i\rightsquigarrow
j}V_j.\eqno{(2)}$$ Here we view $\overline{\tau}$ as an element of
$A$; i. e. $\overline{\tau}(m) = \overline{\tau}\cdot m$.

\begin{laemma} The map $\Phi_{(M,f)}$ is a homomorphism of $A$-modules.
\end{laemma}

\begin{proof} First of all, we recall how $A$ acts on $J$.
Write $J$ as
$$J = \bigoplus_{i\in Q_0}(I_i)\otimes_{\Bbbk}V_i\cong\bigoplus_{i\in Q_0}
\underbrace{I_i\oplus\ldots\oplus I_i}_{\mbox{$\zeta_i$ times}}.$$
For a path $\overline{\sigma}\in\Xi$ denote by $\overline{\tau}^*$
the linear function defined by $\overline{\sigma}^*(\sigma') =
\delta_{\sigma\sigma'}$, for each $\sigma'\in\Xi$, where
$\delta_{\sigma\sigma'}$ stands for the Kronecker delta. As a
$\Bbbk$-linear space $I_i$ has basis
$\overline{\tau}^*_{i1},\ldots,\overline{\tau}^*_{ir(i)}$ with
$\overline{\tau}_{ij}$ being all the paths ending at $i$. To the
basis elements of the $p$-th copy of $I_i$ attach the index $(p)$;
thus, $I_i\otimes V_i = \lspan\left\{\overline{\tau}^{*(p)}_{ij}\mid
p = 1,\ldots,\zeta_i; j = 1,\ldots, r(i) \right\}$. For the image of
an arrow $a$ of $Q$ we have $\overline{a}\cdot\overline{\tau
a}^{*(p)} = \overline{\tau}^{*(p)}$; and
$\overline{a}\cdot\overline{\lambda}^{*(p)} = 0$ in case if $\lambda = e_i$ or $\lambda
= \lambda'b$, where $b$ is an arrow different from $a$. Now pass to
the isomorphism (2). If $a: i\rightarrow k$, then $\overline{a}$
acts on the summand of the right hand side of (2) corresponding to a
path $\overline{\tau}: i\rightsquigarrow j$ (note that the element
$\overline{\tau}^*$ is in $(I_j)_i$) as follows
\begin{align*} &1)\quad \overline{a}: V_j^{(\overline{\tau})}
\xrightarrow{\idd} V_j^{(\overline{\lambda})}, \mbox{ if there is a
path $\lambda$ such that
$\tau = \lambda a$},\\
 &2)\quad \overline{a}:
V_j^{(\overline{\tau})}\rightarrow 0,\mbox{ otherwise}.
\end{align*}
In the first case $\overline{\lambda}: k\rightarrow j$, so that the
image lies in $J_k$. We now check the $A$-invariance of
$\Phi_{M,f}$. Let as before $a: i \rightarrow k$. Then
$$\overline{a}\cdot\varphi_l(m) = \overline{a}(\bigoplus_{\Xi\ni\overline{\tau}: l\rightsquigarrow j}
f_j\overline{\tau}(m)).$$ Now use the above alternative: if
$\overline{\tau} = \overline{\lambda}\overline{a}$ for a path
$\overline{\lambda}: k\rightarrow j$, then the corresponding
component $\varphi_l(m)$ is mapped to the summand
$V^{(\overline{\lambda})}_{j}\subseteq J_k$ without being changed,
and otherwise vanishes. On the other hand,
$$\varphi_k(\overline{a}\cdot m) = \bigoplus_{\Xi\ni\overline{\lambda}: k\rightsquigarrow j}
f_j\overline{\lambda}(\overline{a}m) =
\bigoplus_{\Xi\ni\overline{\lambda}: k\rightsquigarrow j}
f_j(\overline{\lambda}\overline{a})(m) =
\bigoplus_{\begin{smallmatrix}\Xi\ni\overline{\tau}:
l\rightsquigarrow j,\\
\mbox{{\tiny such that $\exists\overline{\lambda}: \overline{\tau} =
\overline\lambda\overline{a}$}}\end{smallmatrix}}
f_j\overline{\tau}(m),$$ which coincides with the above description
of $\overline{a}\cdot\varphi_l(m)$.
\end{proof}

\begin{laemma} The subspace $\ker\Phi_{(M,f)} = \oplus_{i\in
Q_0}\ker\varphi_i$ is the maximal $A$-submodule of $M$ contained in
$\ker{f}$.
\end{laemma}

\begin{proof} It follows from Lemma 2 that $\ker\Phi_{(M,f)}$ is an
$A$-submodule of $M$. 
Now let $U$ be an $A$-submodule of $M$ contained in $\ker{f}$. For
each $\Xi\ni\overline{\tau}:i\rightsquigarrow j$ we then have
$\overline{\tau}U_i = \overline{\tau}e_iU = \overline{\tau} U =
e_j\overline{\tau}U \subseteq U_j$. This implies that
$f_j(\overline{\tau}\cdot x) = 0$, for all $x\in U, j\in Q_0,
\overline{\tau}\in\Xi$, i. e. $U\subseteq\ker\Phi_{(M,f)}$.
\end{proof}

\begin{corol} The map $\Phi_{(M,f)}:M\rightarrow J$ is injective if and only if
the pair $(M,f)$ is stable.
\end{corol}

Now introduce the notion of the {\it Grassmannian of submodules}.
Let $N = \bigoplus_{i\in Q_0}N_i$ be an $A$-module; the Grassmannian
of $A$-submodules of $N$ with dimension vector $\gamma$ is a set
$\Grass^{A}_{\gamma}(N)$ of all $Q_0$-graded subspaces $U =
\bigoplus_{i\in Q_0}U_i\subseteq N$ with $\dim{U_i} = \gamma_i$ that
are $A$-submodules. Note that $\Grass^{A}_{\gamma}(N)$ is a closed
variety in the product of classical Grassmannians
$\Grass_{\gamma}(N) = \prod_{i\in Q_0}\Grass_{\gamma_i}(N_i)$
defined by $\overline{\tau}(U_i)\subseteq U_j$, for all
$\Xi\ni\overline{\tau}:i\rightsquigarrow j$.

The following result is a generalization of~\cite[Prop.
3.9]{Reineke}

\begin{theor} The moduli space $\mathcal{M}^s(A,\alpha,\zeta)$ is
isomorphic to the Grassmannian of submodules $\Grass^A_{\alpha}(J)$.
\end{theor}

\begin{proof} Denote by $\IHom_{\alpha}(V)$ the set of all injective graded
vector space homomorphisms from a space with dimension vector
$\alpha$ to $V$. It is easy to see that $\Grass_{\alpha}(J)$ is a
quotient of $\prod_{i\in Q_0}\IHom_{\Bbbk}(M_i, J_i)$ by the natural
action of $GL(\alpha)$ ($GL(\alpha_i)$ acts on $M_i$ by base
change). Denote the inverse image of
$\Grass^A_{\alpha}(J)\subseteq\Grass_{\alpha}(J)$ by
$\IHom^A_{\alpha}(J)$. Now to prove the theorem it is sufficient to
construct a $GL(\alpha)$-invariant isomorphism
$\Phi:\repq(A,\alpha,\zeta)\xrightarrow{\sim} \IHom^A_{\alpha}(J)$.
We set $\Phi: (M,f)\mapsto\Phi_{(M,f)}$. This map is
$GL(\alpha)$-invariant. Indeed, let $g = (g_i)_{i\in Q_0}\in
GL(\alpha)$. Then for $m\in M$ we have $\Phi(g\cdot(M,f))_i(m_i) =
(\Phi_{g\cdot (M,f)})_i(m) = \bigoplus_{\Xi\ni\overline{\tau}:
i\rightsquigarrow j} f_j(g_j\overline{\tau}m) =
(g\cdot\Phi_{(M,f)})(m)$.

From Corollary 4 it follows that
$\Phi(\repq^s(A,\alpha,\zeta))\subseteq\IHom^A_{\alpha}(J)$. Let us
show that given an injection $\{F: M\hookrightarrow
J\}\in\IHom^A_{\alpha}(J)$, the pair $(M,f)$ may be recovered. But
$f$ is obtained as a composition
$$f:
M_i\xrightarrow{\varphi_i}\bigoplus_{\Xi\ni\overline{\tau}:i\rightsquigarrow
j}V_j\rightarrow V_i,$$ where the last map is a projection on a
summand associated to $\overline{\tau} = e_i$. As for the module
$M$, the following lemma gives the possibility to recover it.

\begin{laemma} Let $K$ be a finite dimensional submodule in $J$ and $\repq^s(A,
\alpha, K)\subseteq\repq^s(A,\alpha)$ be the inverse image of
$\IHom^A_{\alpha}(K)$ under $\Phi$. Then there is a morphism
$s:\IHom^A_{\alpha}(K)\rightarrow\repq^{s}(A,\alpha, K)$ such that
$s\circ\Phi = \idd_{\repq^s(A,\alpha, K)}$.
\end{laemma}

\begin{proof} To recover $M\in\repq^s(A, \alpha, K)$ means to define the action of $A$
on the vector space $\bigoplus_{i\in Q_0}M_i$, i.e. of the elements
$A\ni\overline{a}: i\rightarrow j$ for $a\in Q_1$. Furthermore, for
each of them we have the following commutative diagram:
$$ \xymatrix{
K_i\ar[r]^{K_a} & K_j\\
M_i\ar@{^{(}->}[u]^{f_i}\ar@{-->}[r]^{M_a} &
M_j\ar@{^{(}->}[u]_{f_j}}, $$ where $M_a = \overline{a}|_{M_i}$.
This is because for each pair $(M, f)$ the map $\Phi_{(M,f)}$ is an
$A$-homomorphism. So, each $M_a$ satisfies the equations $\psi_af_i
= f_jM_a$, that may be considered as a matrix equation. Analogously,
$f_j$ may be regarded as a matrix of dimension
$\dim{K_j}\times\dim{M_j}$. Its rank is maximal and equals
$\dim{M_j}$, since $f_j$ are injections. Therefore,
$\IHom^A_{\alpha}(K)$ may be covered by open subsets, where various
minors of the matrix of $f_j$ do not vanish, and $M_a$ are recovered
from matrix elements of $f_i$, $f_j$ and $f_a$ using Kramer's
Theorem.
\end{proof}

Together with the above described way of recovering $f$ the morphism
$s:\IHom^A_{\alpha}(J)\rightarrow\repq^s(A,\alpha)$ gives a morphism
that is inverse to $\Phi$. Consequently, $\Phi$ is an isomorphism
and, being $GL(\alpha)$-invariant, it descends to quotients implying
that $\Phi/\!\!/GL(\alpha):
\mathcal{M}^s(A,\alpha,\zeta)\cong\repq^s(A,\alpha,\zeta)/\!\!/GL(\alpha)\xrightarrow{\sim}
\IHom^A_{\alpha}(J)/\!\!/GL(\alpha)\cong\Grass^A_{\alpha}(J)$.
Theorem 3 is proved.
\end{proof}

Theorem 3 describes a variety that may serve as a substitude of a
moduli space of $A$-modules with dimension vector $\alpha$ whenever
$\repq^s(A, \alpha, \zeta)$ is nonempty. So, it is important to have
a criterion of existence of a stable pair. For quivers M.~Reineke
proved that $\repq^s(A, \alpha, \zeta)\ne\varnothing$ if and only if
$\zeta_i\geqslant(\mathbf{i}, \alpha_i)_Q$, for all $i\in Q_0$,
where $\mathbf{i}$ stands for the vector with all coordinates zero
except $1$ on the $i$-th place, and $(\cdot,\cdot)_Q$ is the Euler
form, i.e. $(\mathbf{i}, \mathbf{j})_Q = \delta_{ij} - (\mbox{number
of arrows from $i$ to $j$})$ , see~\cite[Prop. 4.3]{Reineke}. For
arbitrary finite dimensional algebras we do not have such a result.
However, we can state a weaker proposition. Recall that a socle of
an $A$-module $M$ is the sum of all its simple submodules.

\begin{utver} For an $A$-module $M$ with dimension vector $\alpha$
there is a map $f:M\rightarrow V$ making the pair $(M, f)\in\repq(A,
\alpha, \zeta)$ stable if and only if the socle $\mathrm{soc}\,M$
may be embedded in $\mathrm{soc}\,J$.
\end{utver}

\begin{proof} The ``only if'' part is trivial. We now prove sufficiency. First, note that $(\mathrm{soc}\,J)_i =
V_i^{e_i}$. Consider a decomposition $M = \mathrm{soc}\,M\oplus W$
of $M$ as a vector space. Let $f$ be the composition of the
projection of $M$ onto $\mathrm{soc}\,M$ along $W$ and an inclusion
$\mathrm{soc}\,M\hookrightarrow\mathrm{soc}\,J = \bigoplus_{i\in
Q_0}V_i$. Then $\ker{f} = W$. Now, if $W$ contains a submodule $N$,
then $W\supseteq\mathrm{soc}\,N$, which is a contradiction since
$\mathrm{soc}\,N\subseteq\mathrm{soc}\,M$.
\end{proof}

Taking the dimensions of $(\mathrm{soc}\,J)_i$ and
$(\mathrm{soc}\,M)_i$ his condition may be reformulated as: for
$M\in\repq(Q,\alpha)$ there is a map $f:M\rightarrow V$ making the
pair $(M, f)\in\repq(A, \alpha, \zeta)$ stable if and only if
$\zeta_i\geqslant\dim{(\mathrm{soc}\,M)_i}$. Note that
$(\mathrm{soc}\,M)_i$ is the multiplicity in $\mathrm{soc}\,M$ of
the simple module $S(i) = \Bbbk e_i$ corresponding to the $i$-th
vertex. For quivers this obviously coincides with~\cite[Lemma
4.1]{Reineke}. When using this criterion it is also convenient to
have in mind that $(\mathrm{soc}\,M)_i = \cap_{\overline{a}:
i\rightarrow j}\ker M_{\overline{a}}$.

\medskip

Now, we need a way of determining, for a point of the classical
Grassmannian $\Grass_{\alpha}(J) = \prod_{i\in
Q_0}\Grass_{\alpha_i}(J_i)$, whether it lies in the Grassmanian of
submodules. To formulate the next proposition, we should recall that
the summands $V_s$ of $J_i$ are indexed by paths in $A$.
Furthermore, since an arrow $\overline{a}: i\rightarrow s$ induces,
for each $V^{\overline{\tau}}_k\subseteq J_i$ with $\overline{\tau}
= \overline{\tau}'\overline{a}: i\rightsquigarrow k$, an isomorphism
$V^{\overline{\tau}}_k\xrightarrow{\sim} V^{\overline{\tau}'}_k$,
there is an injection $\dag_{\overline{a}} :
\widetilde{J}^{(\overline{a})}_s\hookrightarrow J_i$, where
$\widetilde{J}^{(\overline{a})}_s$ is a sum of all
$V^{(\overline{\sigma})}_k\subseteq J_s$ such that
$\overline{\sigma}\overline{a}\ne 0$. This injection acts as a
simple index change, and moreover, $\overline{a}\dag_{\overline{a}}
= \idd_{\widetilde{J}^{(\overline{a})}_s}$. Note that for hereditary
algebras, i.e. for quivers with no relations, $\widetilde{J}$
coincides with $J$, so $J_s$ is embedded in $J_i$ whenever there is
an arrow $i\rightarrow s$.

\begin{utver}  A point
$(U_i\subseteq\bigoplus_{\Xi\ni\overline{\tau}i\rightsquigarrow
j}V_j)_{i\in Q_0}\in \prod_{i\in
Q_0}\Grass_{\alpha_i}(\bigoplus_{\Xi\ni\overline{\tau}:i\rightsquigarrow
j}V_j)$ lies in the image of $\Phi/\!\!/GL(\alpha)$ if and only if
$$U_i\subseteq
V_i^{(e_i)}\oplus\dag_{\overline{a}}(U_j\cap
\widetilde{J}^{(\overline{a})}_j),\quad\forall i\in Q_0,\
\overline{a}:i\rightarrow j.\eqno (3)$$
\end{utver}

\begin{proof} As it was shown above, $(U_i\subseteq\bigoplus_{\Xi\ni\overline{\tau}i\rightsquigarrow
j}V_j)_{i\in Q_0}\in \prod_{i\in
Q_0}\Grass_{\alpha_i}(\bigoplus_{\Xi\ni\overline{\tau}:i\rightsquigarrow
j}V_j)$ belongs to the image of $\Phi/\!\!/GL(\alpha)$ if and only
if it is an $A$-submodule of $J$, which means that
$\overline{\tau}(U_i)\subseteq U_j,\,
\forall\Xi\ni\overline{\tau}:i\rightarrow j$. It is straightforward
to check that these conditions are equivalent to (3).
\end{proof}

Denote by $\Bbbk{Q}^{(1)}$ an ideal in the path algebra generated by
all oriented cycles in  $Q$. We also use the notation
$\Bbbk{Q}^{(t)} := (\Bbbk{Q}^{(1)})^t$.

Let $G$ be a reductive group acting on a vector space $X$. Recall,
that the {\it null cone} of this action is the set $\mathcal{N} =
\{x\in X\mid f(x) = 0,\mbox{ for all nonconstant homogeneous
}G\mbox{-invariant functions on X}\}$. By Hilbert-Mumford's
Criterion~\cite[Theorem III.2.4]{Kraft} this is equivalent, for
$x\in X$, to the existence of a one-parameter subgroup
$\lambda:\Bbbk\backslash\{0\}\rightarrow G$ with $\lim_{t\rightarrow
0}\lambda(t)x = 0$.

One of the possible applications of the above construction is the
study of the null cone of $\repq(Q,\alpha,\zeta)$ in case when $Q$
is a quiver with oriented cycles. This is made possible by the
following proposition.

\begin{utver} Let $Q$ be a quiver with oriented cycles.
Then the null cone of the variety $\repq(Q,\alpha,\zeta)$ is
isomorphic to $\repq(\Bbbk Q/(\Bbbk Q)^{(\max_i{\alpha_i})}, \alpha,
\zeta)$.
\end{utver}

\begin{proof} First of all, recall what is the
categorical quotient of a framed representation space. Since there
is a $GL(\alpha)$-invariant isomorphism
$\psi:\repq(Q,\alpha,\zeta)\cong\repq(\widetilde{Q},\widetilde{\alpha})$
for a quiver $\widetilde{Q}$ and an extended dimension vector
$\widetilde{\alpha}$ (see the proof of Proposition 1), we have
$\repq(Q,\alpha,\zeta)/\!\!/GL(\alpha) \cong
\repq(\widetilde{Q},\widetilde{\alpha})/\!\!/GL(\alpha)$ or,
equivalently,
$\psi^*:\Bbbk[\repq(\widetilde{Q},\widetilde{\alpha})]^{GL(\alpha)}\xrightarrow{\sim}
\Bbbk[\repq(Q,\alpha,\zeta)]^{GL(\alpha)}$. But the algebra
$\Bbbk[\repq(\widetilde{Q},\widetilde{\alpha})]^{GL(\alpha)}$ is
generated by the traces of oriented cycles in $\widetilde{Q}$. On
the other hand, oriented cycles in $\widetilde{Q}$ are those of $Q$,
and hence the algebra $\Bbbk[\repq(Q,\alpha,\zeta)]^{GL(\alpha)}$ is
also generated by traces of the oriented cycles in $Q$.

Thus the null cone consists of pairs $(M,f)$, where $M$ are such
representations on which all the oriented cycles in $Q$ act as
nilpotent operators. Since $A^k = 0$ for a nilpotent operator $A$ in
a $k$-dimensional space, all the oriented cycles as operators in $M$
vanish in the $(\max_i{\alpha_i})$-th power. But we state that the
stronger property holds: that the product of any $\max_i{\alpha_i}$
oriented cycles is zero. If $\sigma_1$ and $\sigma_2$ are two such
cycles, then by our conditions $\sigma_1\sigma_2$ and
$\sigma_2\sigma_1$ are also nilpotent, as well as any their product.
Moreover, all of them vanish in the same power. This implies that
the commutator $[\sigma_1,\sigma_2]$ is also nilpotent. Indeed, any
its power equals a sum of products of $\sigma_i$. Since these
products are nilpotent, their traces are all zero, providing that
the trace of any power of the operator $[\sigma_1,\sigma_2]$ is
zero. Therefore, the commutator $[\sigma_1,\sigma_2]$ is nilpotent.

Now, using Engel's Theorem, we can conclude that there is a basis in
the representation space, in which the matrices of both $\sigma_1$
and $\sigma_2$ are upper-triagonal with zeroes on the diagonal.
Applying induction, we get a basis, in which the matrices of all
oriented cycles starting in a given vertex are upper-niltriagonal.
As a product of oriented cycles starting in different vertices is
zero by definition, the product of any $\max_i{\alpha_i}$ such
operators is zero.
\end{proof}

We finish this section with a series of examples.

{\bf Example 1}. Let $A$ be the algebra with quiver $Q : 1
\xrightarrow{a_1} 2\xrightarrow{a_2} 3$ and relations $\rho =
\{a_2a_1\}$. So, $A$ is a Nakayama algebra with admissible sequence
$(1,1)$; see~\cite[Section IV.2]{A-R}. Consider two dimension
vectors $\alpha$ and $\zeta$ and a vector space $V_i$ for each
$\zeta_i$. Then, as we have proved, $\mathcal{M}^s(A, \alpha,
\zeta)\cong\Grass^{A}_{\alpha}(J)$, where $J_1 = V_1\oplus V_2$,
$J_2 = V_2\oplus V_3$, and $J_3 = V_3$. Using Proposition 3 we
obtain that $\mathcal{M}^s(A, \alpha,
\zeta)\cong\left\{(U_i\subseteq J_i)_{i=1}^3\mid U_3\subseteq V_3,
U_2\subseteq V_2\oplus U_3, U_1\subseteq V_1\oplus(U_2\cap
V_2)\right\}$. Note that in this example we do not need to write the
index $(e_i)$ over $V_i$, since there is only one summand $V_i$ in
$J_i$.

{\bf Example 2}. Consider the algebra $A = \Bbbk[x]/(x)^n$. It is
also a Nakayama algebra with admissible sequence (n), but now its
quiver contains a loop. In fact it is a Jordan quiver with a single
vertex and a single loop $a$, and the only relation is $a^n$. Let
$\alpha = (m)$, $\zeta = (k)$, $V$ be a $k$-dimensional vector
space. Then $J = J_1 = V^{(e)}\oplus V^{(a)}\oplus\ldots\oplus
V^{(a^{n-1})}$. Observe that $\overline{a}$ acts on $J$ as follows:
$$\begin{CD} J @. \phantom{U}{:}\phantom{U} @. V^{(e)} @. \phantom{U}\oplus\phantom{U}
 @. V^{(a)} @. \phantom{U}\oplus\phantom{U}
 @. V^{(a^2)} @. \phantom{U}\oplus\phantom{U} @. \ldots @. \phantom{U}\oplus\phantom{U} @. V_{(a^{n-1})}\\
@V{\overline{a}}VV @. @V{0}VV @. @V{\idd}VV @. @V{\idd}VV @. @. @. @V{\idd}VV @. @. @. @.\\
 J @. \phantom{U}{:}\phantom{U} @. 0 @. \phantom{U}\oplus\phantom{U} @.
V^{(e)} @. \phantom{U}\oplus\phantom{U}
 @. V^{(a)} @. \phantom{U}\oplus\phantom{U}
@. \ldots @. \phantom{U}\oplus\phantom{U} @. V^{(a^{n - 2})} @.
\phantom{U}\oplus\phantom{U} @. V^{(a^{n - 1})}
\end{CD}.$$
So, $\widetilde{J} = V^{(e)}\oplus V^{(a)}\oplus\ldots\oplus
V^{(a^{n - 2})}$, and $\dag_{\overline{a}}$ maps $V^{(a^t)}$, $t =
0,\ldots, n - 2$, to $V^{(a^{t + 1})}$ isomorphically.

It is easy to check that $J$ can be viewed as $\Bbbk[x]/(x^n)\otimes
V$ with $\overline{a}$ acting as $\frac{d}{dt}$. Then
$\widetilde{J}$ becomes $\Bbbk[x]/(x^{n-1})\otimes V$,
$\dag_{\overline{a}}: f\otimes v\mapsto xf\otimes v$, and the moduli
space is isomorphic to the Grassmannian of $m$-dimensional
$\frac{d}{dt}$-invariant subspaces in $J$. This interpretation is in
fact rather fruitful and will be further explored in next sections.

{\bf Example 3}. Let $A$ be the algebra with quiver
$$
\xymatrix{
  && 2\ar[dr]^{a_2} & \\
 Q\ : & 1 \ar[ur]^{a_1}\ar[dr]_{b_1} & & 4\\
&& 3\ar[ur]_{b_2} & } $$ and a single relation $a_2a_1 - b_2b_1$.
Then $A$ is a $\Bbbk$-linear span of the path images
$\overline{e}_i$, $i = 1,\ldots,4$, $\overline{a}_i$,
$\overline{b}_i$ and $\overline{a}_2\overline{a}_1$. The components
of $J$ are $J_1 = V_1\oplus V_2 \oplus V_3 \oplus V_4$, $J_2 =
V_2\oplus V_4$, $J_3 = V_3\oplus V_4$ and $J_4 = V_4$ with the arrow
images acting as $\overline{a}_1 = 0\oplus\idd\oplus 0\oplus\idd$,
$\overline{b}_1 = 0\oplus 0\oplus\idd\oplus\idd$, $\overline{a}_2 =
0\oplus\idd$, $\overline{a}_4 = 0\oplus\idd$. Hence $\widetilde{J}_i
= J_i$, for all $i$, all $\dag$'s are natural inclusions. That is
why we venture to omit them as well as the indexes $(e_i)$ in the
final formula. So, $\mathcal{M}^s(A, \alpha,
\zeta)\cong\left\{(U_i\subseteq J_i)_{i=1}^4\mid U_1\subseteq
V_1\oplus U_2, U_1\subseteq V_1\oplus U_3, U_2\subseteq V_2\oplus
U_4, U_3\subseteq V_3\oplus U_4\right\}$.

\section{Framed moduli spaces for quiver $A_{n - 1}^{(1)}$}

In three remaining sections we work over the ground field
$\Bbbk = \mathbb{C}$.

We will see, that in attempts to apply Reineke's construction to
quivers with oriented cycles we have to tackle with Grassmannians in
infinite dimensional spaces, that do not carry an obvious structure
of algebraic variety. However, it is well known that
$\repq^s(Q,\alpha,\zeta)/\!\!/GL(\alpha)$ is projective over
$\repq(Q,\alpha,\zeta)/\!\!/GL(\alpha)$, see~\cite{Kinge}. So, we
can consider the natural projection
$\pi_s:\repq^s(Q,\alpha,\zeta)/\!\!/GL(\alpha)\rightarrow\repq(Q,\alpha,\zeta)/\!\!/GL(\alpha)$
 and investigate its fibers.

In this section our aim is to describe fibers of the projection
$\pi_s$ for a cyclic quiver $Q$ of type $A_{n-1}^{(1)}$:
$$\xymatrix{
1\ar[r]^{a_1} & 2\ar[r]^{a_2}& 3\ar[r]^{a_3}&\ldots\ar[r]^{a_{n -
2}} & (n-1)\ar[dll]^{a_{n - 1}}\\
&& n\ar[ull]^{a_n} &&}$$ These notation for
vertices and arrows of the quiver will be used throughout Sections 4 and 5. 

To each path $\sigma$ in $Q$ we associate a linear function
$\sigma^*\in(\Bbbk Q)^*$, taking values
$$\sigma^*(\sigma') = \begin{cases} 1\mbox{, for $\sigma' = \sigma$},\\
0\mbox{, otherwise}
\end{cases},$$
for a path $\sigma'$. Further, for $\sigma:i\rightsquigarrow j$
denote $V_i^{(\sigma)} := \Bbbk \sigma^*\subseteq(\Bbbk Q)^*$ and
set $\tau_i = a_{i-1}a_{i-2}\ldots a_{i+1}a_i$, the shortest path
starting in $i$. Consider the $\Bbbk Q$-module $J$ with $$(J)_i :=
\prod_{\sigma:i\rightsquigarrow k}V_k^{(\sigma)} = $$ $$ = \prod_{j
= 0}^\infty V_i^{(\tau_i^j)}\times\prod_{j = 0}^\infty V_{i +
1}^{(\tau_{i+1}^ja_i)}\times
\ldots\times\prod_{j = 0}^\infty V_{i - 1}^{(\tau_{i - 1}^ja_{i -
2}\ldots a_{i + 1}a_i)}.\eqno{(4)}$$ Here we identify $\tau_i^0$
with $e_i$. Furthermore, we use the symbol of direct product instead
of direct sum, since we allow our tuples to contain infinitely many
nonzero terms. It should be explained how $a_i\in\Bbbk Q$ acts on
$J$. In order to do this we use the alternative formulated in the
proof of Lemma 2. Considering the restrictions of $a_i$ to the
summands of (4) we have
$$a_i: J_i\supseteq\prod_{j =
0}^\infty V_k^{(\tau_k^ja_{k - 1}\ldots a_{i +
1}a_i)}\xrightarrow{\quad\sim\quad}\prod_{j = 0}^\infty
V_k^{(\tau_{k}^ja_{k - 1}\ldots a_{i + 1})}\subseteq J_{i + 1},$$
for $k\ne j$, where $a_i$ acts as componentwise isomorphism
$V_k^{(\tau_k^ja_{k - 1}\ldots a_{i +
1}a_i)}\xrightarrow{\idd}V_k^{(\tau_{k}a_{k - 1}\ldots a_{k + 1})}$.
In cases if $k = i$ a shift $J_i\supseteq
V_i^{(\tau_i^j)}\xrightarrow{\idd} V_k^{(\tau_{i}^{j - 1}a_{i -
1}\ldots a_{i + 1})}\subseteq J_{i + 1}$ takes place and, in
particular, $a_i: J_i\supseteq V_i^{(e_i)}\rightarrow 0$.

For a pair $(M, f)\in\repq(Q,\alpha,\zeta)$ define as before the map
$\Phi_{(M,f)} = (\varphi_i)_{i = 1}^n: M\rightarrow J$ by
$\varphi_i(m) = (f_j(\tau m))_{\tau: i\rightsquigarrow j}$. We mean
here that the component of $\varphi(m)$ that corresponds to
$V_k^{(\tau)}\subseteq J_i$ equals $f(\tau m)$.

{\bf Example 4}. If $n = 1$, then $J = J_1 = V^{(e_1)}\times
V^{(a)}\times V^{(a^2)}\times\ldots\times V^{(a^k)}\times\ldots$ ($a
= a_1: i\rightarrow i$ is the only arrow in the quiver) and the map
$\Phi_{(M,f)}$ is of form $\varphi(m) = (f(m), f(am),
f(a^2m),\ldots,f(a^km),\ldots)$. If $\zeta_1 = \dim{V_1} = 1$, then
$J$ is isomorphic to the ring $\Bbbk[[t]]$ of formal power series
with coefficients in $\Bbbk$. The isomorphism is realized by $T:
(c_i)_{i = 0}^\infty\mapsto\sum_{i = 0}^\infty\frac{c_i}{i!}t^i$,
where the arrow $a$ acts as $\frac{d}{dt}$ in $\Bbbk[[t]]$. Indeed,
$\frac{d}{dt}T((c_i)_{i = 0}^\infty) = \frac{d}{dt}(\sum_{i =
0}^\infty)\frac{c_i}{i!}t^i = \sum_{i = 0}^\infty\frac{c_{i +
1}}{i!}t^i = T((c_{i + 1})_{i = 0}^\infty) = T(a\cdot(c_i)_{i =
0}^{\infty})$, hence $T$ is really an isomorphism of $\Bbbk
Q$-modules. We may further improve $\Phi_{(M,f)}$ by setting
$\Phi^!_{(M,f)} := T\circ\Phi_{(M,f)}: m\mapsto\sum_{i =
0}^\infty\frac{f_1(a_1^im)}{i!}t^i = f_1(\sum_{i =
0}^\infty\frac{1}{i!}a_1^it^im) = f_1(\exp(a_1t)m)$.

\medskip

Inspired by this example, we'll try to generalize this
interpretation for arbitrary $n$ and $\zeta$. Namely, consider the
space $J^! = \bigoplus_{i = 1}^nJ^!_i$ with $J^!_i =
(\Bbbk[[t]]\otimes V_1)\oplus(\Bbbk[[t]]\otimes
V_2)\oplus\ldots\oplus(\Bbbk[[t]]\otimes V_n)$, for $i =
1,\ldots,n$, with the following $\Bbbk Q$-module structure:
$a_i\cdot(G_1,\ldots,G_n) = (G_1, \ldots, G_{i - 1},\frac{d}{dt}G_i,
G_{i + 1},\ldots, G_{n})$ (elements of $\Bbbk[[t]]\otimes V_i$ can
be regarded as vectors $G = (g_1,\ldots,g_{\zeta_i})$,
$g_j\in\Bbbk[[t]]$; in this situation $\frac{d}{dt}$ is a
componentwise formal differentiation). Further, define a map $T =
(T_i)_{i = 1}^n: J\rightarrow J^!$ by
$$T_i: J_i = \prod_{j =
0}^\infty V_i^{(\tau_i^j)}\times\prod_{j = 0}^\infty V_{i +
1}^{(\tau_{i+1}^ja_i)}\times
\ldots\times\prod_{j = 0}^\infty V_{i -
1}^{(\tau_{i - 1}^ja_{i - 2}\ldots a_{i + 1}a_i)}\longrightarrow$$
$$\longrightarrow J^!_i = (\Bbbk[[t]]\otimes
V_i)\oplus(\Bbbk[[t]]\otimes V_{i + 1})\oplus(\Bbbk[[t]]\otimes V_{i
+ 2})\oplus\ldots\oplus(\Bbbk[[t]]\otimes V_{i - 1});$$ this map is
componentwise and its components $T_{ik}: \prod_{j = 0}^\infty
V_k^{(\tau_{k}^ja_{k - 1}\ldots a_{i +
1}a_i)}\rightarrow(\Bbbk[[t]]\otimes V_k)$ are
$$T_{ij}(\overline{c}_0, \overline{c}_1,\ldots,\overline{c}_k,\ldots)
= \sum_{k = 0}^\infty\frac{1}{k!}t^k\overline{c}_k.$$

\begin{laemma} The map $T$ is a $\Bbbk Q$-isomorphism.
\end{laemma}

\begin{proof} We need to show that for every $i$, $k$ there is a
commutative diagram of restrictions:
$$\begin{CD} J_i @. \phantom{U}\subseteq\phantom{U} @. \prod_{j =
0}^\infty V_k^{(\tau_k^ja_{k - 1}\ldots a_{i + 1}a_i)} @>{a_i}>>
\prod_{j = 0}^\infty
V_k^{(\tau_{k}^ja_{k - 1}\ldots a_{k + 1})} @. \phantom{U}\subseteq\phantom{U} @. J_{i + 1}\\
@. @. @V{T_{ik}}VV @V{T_{i + 1, k}}VV @. @.\\
 J^!_i @. \phantom{U}\subseteq\phantom{U} @. (\Bbbk[[t]]\otimes
V_k) @>{a_i}>> (\Bbbk[[t]]\otimes V_k) @.
\phantom{U}\subseteq\phantom{U} @. J^!_{i+ 1}
\end{CD}.$$
Observe that for $k\ne i$ the restriction of $a_i$ to the
subspaces considered is the identity operator in both rows, hence
the diagram is commutative. Otherwise, if $k = i$, then $a_i\cdot
T_{ii}(\overline{c}_0, \overline{c}_1,\ldots,\overline{c}_k,\ldots)
= \frac{d}{dt}\left(\sum_{k =
0}^\infty\frac{1}{i!}t^k\overline{c_i}\right) = \sum_{k =
0}^\infty\frac{1}{i!}t^k\overline{c_{i + 1}} = T_{i + 1,
i}(\overline{c}_1, \overline{c}_2,\ldots,\overline{c}_{k +
1},\ldots) = T_{i, i + 1}(a\cdot(\overline{c}_0,
\overline{c}_1,\ldots,\overline{c}_k,\ldots))$. This finishes the
proof.
\end{proof}

Now, just as we did in the last example we turn to maps
$\Phi^!_{(M,f)} = (\varphi^!_i)^n_{i = 1} := T\circ\Phi_{M,f}$.
Having in mind the decomposition (4), it is not hard to see that
\begin{align*}
\varphi^!_{i}(m) =\ &(\sum_{k =
0}^\infty\frac{1}{k!}t^kf_i(\tau_i^km), \sum_{k =
0}^\infty\frac{1}{k!}t^kf_{i + 1}(\tau_{i + 1}^ka_im),\sum_{k =
0}^\infty\frac{1}{k!}t^kf_{i + 2}(\tau_{i + 2}^ka_{i +
1}a_im), \\
 &,\ldots,\sum_{k = 0}^\infty\frac{1}{k!}t^kf_{i - 1}(\tau_{i
- 1}^ka_{i-2}\ldots a_{i + 1}a_im)) =
\\
\!\!\!= &(f_i(\exp(\tau_it)m), f_{i + 1}(\exp(\tau_{i + 1}t)a_im),
f_{i + 2}(\exp(\tau_{i + 2}t)a_{i + 1}a_im),
\ldots,\\
&\phantom{AA}f_{i - 1}(\exp(\tau_{i - 1}t)a_{i-2}\ldots a_{i +
1}a_im)).
\end{align*}

Now it will be shown that some of the useful properties of its
finite dimensional analog hold for the above constructed map.

\begin{laemma} $\phantom{AAA}$
\begin{itemize}
\item[(i)] The map $\Phi^!_{(M,f)}$ is a $\Bbbk
Q$-homomorphism.
\item[(ii)] The subspace $\ker{\Phi^!_{(M,f)}}$ is the maximal
$\Bbbk Q$-submodule of $M$ contained in $\ker{f}$.
\item[(iii)] A pair $(M,f)$ is stable if and only if
$\Phi^!_{(M,f)}$ is injective.
\end{itemize}
\end{laemma}

\begin{proof} (i) Let $a_i$ be an arrow in $Q$; denote by $\varphi^!_i$
the components $\Phi^!_{(M,f)}$. It is clear that it is sufficient
to show that $(a_i)|_{J^!_i}\circ\phi^!_i = \phi^!_{i +
1}\circ(a_i)|_{M_i}$.

We have $a_i\cdot\varphi^!_{i}(m) = a_i\cdot(f_i(\exp(\tau_it)m),
f_{i + 1}(\exp(\tau_{i + 1}t)a_im),
\ldots,$ $f_{i - 1}(\exp(\tau_{i -
1}t)a_{i-2}\ldots a_{i + 1}a_im)) =
(\frac{d}{dt}f_i(\exp(\tau_it)m), f_{i + 1}(\exp(\tau_{i +
1}t)a_im),
$ $\ldots,\!f_{i - 1}(\exp(\tau_{i - 1}t)a_{i-2}\ldots a_{i +
1}a_im)) \!=\! (\!f_i(\exp(\tau_it)\tau_im),\!f_{i + 1}(\exp(\tau_{i
+ 1}t)a_im),
$ $\ldots,f_{i - 1}(\exp(\tau_{i - 1}t)a_{i-2}\ldots a_{i + 1}a_im))
= (f_i(\exp(\tau_it)a_{i-1}\ldots a_{i + 1}\cdot a_im),$ $f_{i +
1}(\exp(\tau_{i + 1}t)\cdot a_im),
\!\ldots,\!f_{i - 1}(\exp(\tau_{i - 1}t) a_{i-2}\!\ldots\!a_{i +
2}a_{i + 1}\cdot a_im))\!=\!\phi^!_{i + 1}(a_im)$. Thus,
$\Phi^!_{(M,f)}$ is a $\Bbbk Q$-homomorphism.

(ii) It follows from the first part of this lemma that
$\ker{\Phi^!_{(M,f)}}$ is a $\Bbbk Q$-submodule of $M$. Next, let
$N\subseteq M$ be a submodule contained in $\ker{f}$. Then for a
path $\tau$ we have $\tau(N)\subseteq N$ and, therefore,
$\Phi_{(M,f)}(N) = 0$ (it is a straightforward consequence of our
definitions). Consequently, $N\subseteq\ker{\Phi_{(M,f)}} =
\ker{\Phi^!_{(M,f)}}$.

(iii) It is a direct consequence of (ii).
\end{proof}

In the finite dimensional case to each element
$(M,f)\in\repq^s(Q,\alpha,\zeta)$ associated is a submodule of $J^!$ with
dimension vector $\alpha$, that is a point in the Grassmannian of
submodules $\Grass_{\alpha}^{\Bbbk Q}(J^!)$. However, to obtain a
one-to-one correspondence, we need the following simplification.

 Denote by
$\mathcal{A}[t]$ the subspace in $\Bbbk[[t]]$ consisting of power
series converging for all $t$. Since $\Bbbk = \mathbb{C}$, this
conditions means that the series is a Taylor series of a holomorphic
function, and hence the uniqueness theorem implies that
$\mathcal{A}[t]$ is the ring of entire functions
$\mathcal{O}(\mathbb{C})$. As $\varphi^!_i$ are expressed in terms
of exponents, they all lie in $J^{!!}_i := (\mathcal{A}[t]\otimes
V_i)\oplus(\mathcal{A}[t]\otimes V_{i +
1})\oplus(\mathcal{A}[t]\otimes V_{i +
2})\oplus\ldots\oplus(\mathcal{A}[t]\otimes V_{i - 1})$. Set $J^{!!}
:= \bigoplus_{i = 1}^nJ^{!!}_n$. It is not hard to see that $J^{!!}$
is a $\Bbbk Q$-submodule in $J^{!}$, and so our task is now to
describe all $\Bbbk Q$-submodules in $J^{!!}$ with dimension vector
$\alpha$, i.e. the set $\Grass^{\Bbbk Q}_{\alpha}(J^{!!})$ and to
prove that such submodules are in one-to-one correspondence with
stable pairs
$(M,f)\in\repq^s(Q,\alpha,\zeta)$. 

%
%

Let $U\subseteq J^{!!}$ be a $\Bbbk Q$-submodule with dimension
vector $\alpha$. Then for each $i = 1,\ldots,n$ we have an inclusion
$a_i(U_i)\subseteq U_{i + 1}$. This in particular implies that
$\tau_{i}(U_i)\subseteq U_i$. But we know that $\tau_i: J^{!!}_i =
\bigoplus_{j =r 1}^n(\mathcal{A}[t]\otimes V_j)\rightarrow
\bigoplus_{j = 1}^n(\mathcal{A}[t]\otimes V_j) = J^{!!}_{i + 1},\
(G_1,\ldots,G_n)\mapsto(\frac{d}{dt}G_1,\ldots,\frac{d}{dt}G_n)$.
Hence if $U$ is a submodule, then it is preserved by the
componentwise differentiation of elements $G\in U_i$. This means
that $U_i$ are $\frac{d}{dt}$-invariant subspaces in $J^{!!}_i$
(note that $\bigoplus_{j = 1}^n(\mathcal{A}[t]\otimes
V_j)\cong\mathcal{A}[t]\otimes(\bigoplus_{j = 1}^nV_j)$ and so its
elements may be considered as $\sum_{j = 1}^n\zeta_j$-tuples of
functions). Further, let $D^{(i)} = (d^{(i)}_{pq})$ be the matrix of
the restricted operator $\frac{d}{dt}|_{U_i}$ and set $k =
\alpha_i$, $m = \sum_{j = 1}^n\zeta_j$. Fix a base
$\overline{g}_1,\ldots,\overline{g}_{k}$ in $U_i$ (following the
idea expressed in a recent remark we here consider $\overline{g}_j$
as a $m$-tuple of functions: $\overline{g}_j = (g_{jl})_{l = 1}^n$).
Then, for all $j$, we get:
$$\frac{d}{dt}\overline{g}_j = d^{(i)}_{1j}\overline{g}_1 +
d^{(i)}_{2j}\overline{g}_2 + \ldots + d^{(i)}_{kj}\overline{g}_k.$$
It is an easy calculation to check that
$$\left(\overline{g}_1(t),\ldots,\overline{g}_k(t)\right) =
\left(\overline{g}_1(0),\ldots,\overline{g}_k(0)\right)\exp(D^{(i)}t).\eqno{(5)}$$
Now Cayley-Hamilton's Theorem implies that $\chi_i(D^{(i)}) = 0$ in
$U_i$, where $\chi_i$ is a characteristic polynomial of $D^{(i)}$.
Thus each vector function $\overline{g}\in U_i$ satisfies the
differential equation $\chi_i(\frac{d}{dt})\overline{g} = 0$.

We can finally show that the above constructed correspondence
between stable pairs and points in Grassmannian is one-to-one.

\begin{utver} The map $\Phi^{!}: \repq^s(Q,\alpha,\zeta)\rightarrow
\IHom^{\Bbbk Q}_{\alpha}(J^{!!}),\ (M,f)\mapsto\Phi^{!}_{(M,f)}$ is
a bijection.
\end{utver}

\begin{proof} We need to show that having an inclusion $\Phi^{!}\in
\IHom^{\Bbbk Q}_{\alpha}(J^{!!})$ we can uniquely recover a pair
$(M,f)$. First of all note that since, for each $j$,$r$, we have
$\tau_j^ra_{j - 1}\ldots a_i = a_{j - 1}\ldots a_i\tau_i^r$, the map
$\varphi^!_{i}$ may be written as:
$$
\varphi^!_{i}(m) = (f_i(\exp(\tau_it)m), f_{i +
1}(a_i\exp(\tau_{i}t)m),\ldots,$$ $$f_{i - 1}(a_{i-2}\ldots a_{i +
1}a_i\exp(\tau_{i}t)m)).\eqno{(6)}
$$

Let now $U = (U_i)_{i = 1}^n\subseteq J^{!!}$ be a submodule with
dimension vector $\alpha$. Let also $U_i =
\lspan\left\{\overline{g}^{(i)}_1(t),\ldots,\overline{g}^{(i)}_{\alpha_i}(t)\right\}$.

Recall the equality (5), setting $\overline{g}_r(t) =
\overline{g}^{(i)}_r(t)$. The matrices
$\left(\overline{g}^{(i)}_1(t),\ldots,\overline{g}^{(i)}_{\alpha_i}(t)\right)$
and
$\left(\overline{g}^{(i)}_1(0),\ldots,\overline{g}^{(i)}_{\alpha_i}(0)\right)$
may be divided into horizontal blocks of size
$\zeta_j\times\alpha_i$, those blocks corresponding to the natural
projections $P_j: J^{!!}_i\twoheadrightarrow\mathcal{A}[t]\otimes
V_j$. Thus, for each $j$, we obtain
\begin{equation*}\left(P_j\overline{g}^{(i)}_1(t),\ldots,P_j\overline{g}^{(i)}_{\alpha_i}(t)\right)
=
\left(P_j\overline{g}^{(i)}_1(0),\ldots,P_j\overline{g}^{(i)}_{\alpha_i}(0)\right)
\exp{D^{(i)}}t\cdot E_{\alpha_i},\eqno{(7)}
\end{equation*}
where $E_{\alpha_i}$ is the identity matrix of size $\alpha_i\times
\alpha_i$. These equalities can be interpreted as follows: there is
a map $\Psi = (\Psi_{ij})^{n}_{i,j = 1}$ with $\Psi_{ij} =
\left(P_j\overline{g}^{(i)}_1(0),\ldots,P_j\overline{g}^{(i)}_{\alpha_i}(0)\right)
\exp{D^{(i)}}t: M_i\rightarrow\mathcal{A}[t]\otimes V_j\subseteq
J^{!!}_i$, where $M$ is a graded vector space with dimension vector
$\alpha$; one may easily establish that $\Psi$ is bijective. We need
to show that $\Psi = \Phi^{!}_{(M,f)}$ for a map $f: M\rightarrow V$
and a certain $\Bbbk Q$-module structure on $M$.

As for the module structure, it is quite clear: viewing (6) as
formulas defining the natural inclusion $U\hookrightarrow J^{!!}$,
we set $M = U$. Further, $f = (f_i)_{i = 1}^n$ is defined as a tuple
of compositions $f_i: M_i \hookrightarrow J^{!!}_i
\twoheadrightarrow V_i$, where $V_i$ should be regarded as the
component of $J^{!!}_i$ associated to $e_i$ (we mean
$V_i^{(e_i)}\subseteq V_i^{(e_i)}\oplus\bigoplus_{j =
1}^{\infty}V_i^{\tau_i^j} = (\Bbbk[[t]]\otimes V_i)$). Hence,
$$f_i = \Psi_{ii}|_{t = 0}: m\mapsto
\left(P_i\overline{g}^{(i)}_1(0),\ldots,P_i\overline{g}^{(i)}_{\alpha_i}(0)\right)m.$$

Now it is left to prove that (7) defines a map $\Psi$ coinciding
with $\Phi_{(M,f)}$ for the above $M$ and $f$. From (6) it follows,
that it is sufficient to show, that
$$\left(P_j\overline{g}^{(i)}_1(0),\ldots,P_j\overline{g}^{(i)}_{\alpha_i}(0)\right)
\exp{D^{(i)}}t = f_j\cdot a_{j - 1}\ldots a_i\exp(\tau_it),\eqno(8)
$$
for $j\ne i$. But $D^{(i)}$ is just a notation for the matrix of
$\tau_i$, so $\exp(D^{(i)}t)\equiv\exp(\tau_it)$. Canceling this
exponent, we turn (8) into
$$\left(P_j\overline{g}^{(i)}_1(0),\ldots,P_j\overline{g}^{(i)}_{\alpha_i}(0)\right)
 = \left(P_j\overline{g}^{(j)}_1(0),\ldots,P_j\overline{g}^{(j)}_{\alpha_j}(0)\right)
 \cdot a_{j - 1}\ldots a_i.$$

We are going to prove this as a matrix equality, and so we can
immediately write
$$\left(a_{j - 1}\ldots a_i\overline{g}^{(i)}_1(t),\ldots,a_{j - 1}\ldots a_i\overline{g}^{(i)}_{\alpha_i}(t)\right)
=
\left(\overline{g}^{(j)}_1(t),\ldots,\overline{g}^{(j)}_{\alpha_j}(t)\right)
\cdot a_{j - 1}\ldots a_i,$$ where $a_{j - 1}\ldots
a_i\overline{g}^{(i)}_r(t)$ stands for the value of $a_{j - 1}\ldots
a_i\overline{g}^{(i)}_r$ in $t$. Descending to the level of
projections, we obtain
$$\left(P_ja_{j - 1}\ldots a_i\overline{g}^{(i)}_1(t),\ldots,P_ja_{j - 1}\ldots a_i\overline{g}^{(i)}_{\alpha_i}(t)\right)
= $$
$$ = \left(P_j\overline{g}^{(j)}_1(t),\ldots,P_j\overline{g}^{(j)}_{\alpha_j}(t)\right)
\cdot a_{j - 1}\ldots a_i,$$ but $a_r$ acts on $J^{!!}_i$ as
$(\idd,\ldots,\frac{d}{dt},\ldots,\idd)$, where the differentiation
takes place only at the $r$-th position. Consequently, the product
$a_{j - 1}\ldots a_i$ does not change the $j$-th projection of
$\overline{g}^{(i)}$, implying that
$$\left(P_ja_{j - 1}\ldots a_i\overline{g}^{(i)}_1(t),\ldots,P_ja_{j - 1}\ldots a_i\overline{g}^{(i)}_{\alpha_i}(t)\right)
=
\left(P_j\overline{g}^{(i)}_1(t),\ldots,P_j\overline{g}^{(i)}_{\alpha_i}(t)\right),$$
and the required equality follows. The proposition is proved.
\end{proof}

\begin{corol} Points of $\mathcal{M}^s(Q,\alpha,\zeta)$ are in
one-to-one correspondence with points of the Grassmannian of
submodules $\Grass_{\alpha}^{\Bbbk Q}(J^{!!})$.
\end{corol}

This, however, does not give us a desired isomorphism. The reason is
that instead of describing the quotient we rather have proved that
$\Grass_{\alpha}^{\Bbbk Q}(J^{!!})$ is an algebraic variety.

Now consider the standard categorical quotient
$\pi_{GL(\alpha)}:\repq(Q,\alpha,\zeta)\rightarrow\mathcal{M}(Q,
\alpha) = \repq(Q, \alpha,\zeta)/\!\!/GL(\alpha)$. As it was pointed
out before, this quotient parameterizes all possible characteristic
polynomials of the cycles $\tau_1,\ldots,\tau_n$ in $Q$; following
this observation, we will consider $\mathcal{M}(Q, \alpha,\zeta)$ as
embedded in the product
$\Bbbk[x]_{\alpha_1}\times\ldots\times\Bbbk[x]_{\alpha_n}$, where
$\Bbbk[x]_r$ is a space of polynomials of degree nor higher than $r
+ 1$. It is clear that not each tuple of polynomials can be obtained
from a representation. Although we are not going to give an explicit
description of $\mathcal{M}(Q, \alpha,\zeta)$ as a subvariety in
$\Bbbk[x]_{\alpha_1}\times\ldots\times\Bbbk[x]_{\alpha_n}$, we prove
the following useful lemma.

\begin{laemma} Let $\overline{\chi} = (\chi_1,\ldots,\chi_n)$ be an
admissible tuple of polynomials and $\lambda_1,\ldots,\lambda_N$ be
all different roots of $\chi_1,\ldots,\chi_n$. Let also $r_{ij}$ be
the multiplicity of $\lambda_j$ as a root of $\chi_i$. If
$\lambda_j\ne 0$, then $r_{1j} = \ldots = r_{nj}$.
\end{laemma}

\begin{proof} Let $M\in\pi_{GL(\alpha)}^{-1}(\overline{\chi})$.
Consider the Jordan block decomposition of $M_i$ with respect to
$\tau_i$. Then $r_{ij}$ is the dimension of the subspace
$V_i^{\lambda_j} := \left\{m\in M_i\mid\exists q\in\mathbb{N}:
(\tau_i - \lambda_j\idd)^qm = 0\right\}$. Since the restriction of
$\tau_i = a_{i-1}\ldots a_i$ on $V_i^{\lambda_j}$ is non-degenerate,
$a_{p-1}\ldots a_i(M_i^{\lambda_j})\cap\ker{a_p} = 0$, for all $p =
1,\ldots,n$. As $\tau_pa_{p - 1}\ldots a_i = a_{p - 1}\ldots
a_i\tau_i$, all $\tau_p|_{a_{p-1}\ldots a_i(M_i^{\lambda_j})}$ are
conjugate and therefore have the same eigenvalues.
\end{proof}

Recall that there is a natural projection $\pi_s:\mathcal{M}^s(Q,
\alpha,\zeta)\rightarrow\mathcal{M}(Q, \alpha,\zeta)$. We are going
to investigate its fibers $\pi_s^{-1}(\overline{\chi})$, for
$\overline{\chi} = (\chi_1,\ldots,\chi_n)\in\mathcal{M}(Q, \alpha,
\zeta)$.

Fix a tuple $\overline{\chi}$. Let $\lambda_j$ and $r_{ij}$ be as in
Lemma 7. Define $\underline{r}_j = (r_{1j}, \ldots, r_{nj})$, for $j
= 1,\ldots,N$. Consider the submodules $J(\lambda_j,
\underline{r}_j)$ of $J^{!!}$, where $J(\lambda_j,
\underline{r}_j)_i$ is the subspace in
$\mathcal{A}[t]\otimes(\bigoplus_{j = 1}^nV_j)$ generated by all
solutions of the differential equation $(\frac{d}{dt} -
\lambda_j\idd)^{\max_i{r_{ij}}}\cdot\overline{g} = 0$.

\begin{laemma} For $M\!\in\!\pi_s^{-1}(\overline{\chi})$, the
image $\Phi^!_{(M,f)}(M)$ lies in $\bigoplus_{j = 1}^NJ(\lambda_j, \underline{r}_j)$. 
\end{laemma}

\begin{proof} As was shown before,
each vector function in $(\Phi_{(M,f)}(M))_i$ satisfies the equation
$\chi_i(\frac{d}{dt})\cdot\overline{g} = 0$, where $\chi_i$ is a
characteristic polynomial of $\tau_i$ as an operator in $M_i$. If
$\chi_i(t) = (t - \lambda_1)^{r_{i1}}\ldots(t -
\lambda_N)^{r_{iN}}$, then clearly $(\Phi_{(M,f)}(M))_i\subseteq
J(\lambda_1, \underline{r}_1)_i\oplus\ldots\oplus J(\lambda_N,
\underline{r}_N)_i$. This yields the required inclusion.
\end{proof}

Let now $W = \bigoplus_{j = 1}^NJ(\lambda_j, \underline{r}_j)$ and $\Phi^j_{(M, f)}$ be $\Phi^!_{(M, f)}$ followed by the
projection on $J(\lambda_j, \underline{r}_j)$ along $\bigoplus_{p\ne
j}J(\lambda_p, \underline{r}_p)$.

\begin{laemma} A stable pair $(M, f)$ is in $\pi^{-1}_{GL(\alpha)}(\overline{\chi})$
if and only if the dimension vector of each $\Phi^j_{(M,f)}(M)\subseteq
J(\lambda_j, \underline{r}_j)$, $j = 1,\ldots, N$, equals
$\underline{r}_{ij}$.
\end{laemma}

\begin{proof} The ``only if'' part is straightforward. Conversely, if $(M, f)$ is
a stable pair, we have $\Phi^j_{(M,f)}(M)
= \bigoplus_{j = 1}^NW(j)$, where $W(j) = W\cap J(\lambda_j, \underline{r}_j)$, so
that $W(j)_i = W_i^{\lambda_j}$. Since $\Phi^!_{(M,f)}$ is a $\Bbbk
Q$-homomorphism, $\dim{M_i^{\lambda_j}} = \dim{W_i^{\lambda_j}} =
r_{ij}$. Therefore the multiplicity of $\lambda_j$ as an eigenvalue
of $\tau_i$ as an operator on $M_i$ equals $r_{ij}$. This implies
that the characteristic polynomial of $\tau_i$ acting on $M_i$ is
$\prod_{j = 1}^N(t - \lambda_j)^{r_{ij}} = \chi_i(t)$. So, $(M,
f)\in\pi^{-1}_{GL(\alpha)}(\overline{\chi})$.
\end{proof}

This lemma ensures that $\Phi^!:
\repq^s(Q,\alpha,\zeta)\rightarrow\IHom_{\alpha}^{\Bbbk Q}(J^{!!}),\
(M, f)\mapsto\Phi^!_{(M, f)}$ restricted to
$\pi_{GL(\alpha)}^{-1}(\overline{\chi})$ induces a bijection
$\repq^s(Q,\alpha,\zeta)\cap\pi_{GL(\alpha)}^{-1}(\overline{\chi})\rightarrow
\prod_{j = 1}^N\IHom_{\underline{r}_j}^{\Bbbk Q}(J(\lambda_j,
\underline{r}_j))$. As $J(\lambda_j, \underline{r}_j)$ are finite
dimensional, we can use Lemma 4 to prove that
$\pi_s^{-1}(\overline{\chi})$ is isomorphic (this time as an
algebraic variety) to $\prod_{j =
1}^N\Grass_{\underline{r}_j}^{\Bbbk Q}(J(\lambda_j,
\underline{r}_j))$.

Collecting the results obtained we can state the following:

\begin{theor}
\begin{itemize}
\item[(1)]
Points of the quotient $\mathcal{M}^s(Q,\alpha,\zeta)$ are in
one-to-one correspondence with points of the Grassmannian of
submodules $\Grass_{\alpha}^{\Bbbk Q}(J^{!!})$.
\item[(2)] Let $\overline{\chi} = (\chi_1,\ldots,\chi_n)$ be an
admissible tuple of polynomials, $\lambda_1,\ldots,\lambda_N$ be all
different roots of $\chi_1,\ldots,\chi_n$, and $r_{ij}$ be the
multiplicity of $\lambda_j$ as a root of $\chi_i$. Then
$\pi_s^{-1}(\overline{\chi})\cong\prod_{j =
1}^N\Grass_{\underline{r}_j}^{\Bbbk Q}(J(\lambda_j,
\underline{r}_j))$, where $\underline{r}_j = (r_{1j}, \ldots,
r_{nj})$.
\end{itemize}
\end{theor}

\section{An explicit realization of fibers}

Let $Q$ be a Jordan quiver consisting of a single vertex and a
single loop (both this loop and the corresponding operator in a
representation will be denoted by $a$). Set also $\alpha = (m)$,
$\zeta = (q)$. This is the case when our construction becomes as
clear as possible.

It is evident that the standard categorical quotient for the action
$GL(m) : \repq(Q,m)$ is isomorphic to $\mathbb{A}^m$: points of the
quotient are tuples of characteristic polynomial coefficients of the
operator corresponding to the arrow $a$; having this in mind we will
further assume that $\mathcal{M}(Q, m)$ is embedded in $\Bbbk[t]_m$.

We have $J^{!!} = J^{!!}_1 = \mathcal{A}\otimes V_1$, and the map
$\Phi^!_{(M,f)}$ becomes $\varphi(m) = f_1(\exp(at)m)$. Further, for
a subspace $U\in\Grass_m(J^{!!})$ the equivalence
$U\in\TrueIm(\Phi^!)\Leftrightarrow \frac{d}{dt}(U)\subseteq U$
holds. Thus, 
the fiber over
$\chi\in\Bbbk[x]$ is precisely
$\Grass^{\frac{d}{dt}}_{m}(J^{!!}_{\chi})$ (the Grassmannian of
$m$-dimensional $\frac{d}{dt}$-invariant subspaces in
$J^{!!}_{\chi}$).

Imagine $q = 1$. We then have ordinary functions instead of vector
ones; and the dimension of the solution space of the differential
equation $\chi(\frac{d}{dt})g = 0$ equals $m$. So each Grassmannian
$\Grass_{}^{\frac{d}{dt}}(J^{!!}_{\chi})$ consists in this case of a
single subspace. It other words, in each fiber of the projection
$\pi_s: \mathcal{M}^s(Q,m,1)\twoheadrightarrow\mathcal{M}(Q,m,1)$
there is at most one point. Such a result is rather upsetting,
though the situation will be more favorable for $q
> 1$. We may even guarantee that for $q = m$ each $\Bbbk Q$-module $M$ arises as a member of
a stable pair $(M,f)\in\repq^s(Q,m,q)$ (it is, for instance, $(M,
\idd)$).

The next theorem describes the fibre structure of $\pi_s$ for
arbitrary $q$ and $m$.

\begin{theor} Let $m, q$ be positive integers.
\begin{itemize}
\item[(a)] If $\chi(x) = \prod_{i = 1}^s(x - \lambda_i)^{r_i}$ and all
$\lambda_i$ are different, then $\pi_s^{-1}(\chi)\cong\prod_{i =
1}^s\Grass_{r_i}^{\frac{d}{dt}}(J(\lambda_i, r_i))$.
\item[(b)] If $m = 1$, then $\Grass_m^{\frac{d}{dt}}(J(\lambda, m))\cong\mathbb{P}^{q -
1}$.
\item[(c)] If $m > 1$ and $\chi(x) = (x - \lambda)^m$, then $\Grass_m^{\frac{d}{dt}}(J(\lambda, m))$
is isomorphic to the subvariety in $\mathbb{P}(\bigwedge^mJ(\lambda,
m))$ given by the following equations:
\begin{itemize}
\item[(i)] The Plucker equations;
\item[(ii)] $p_{i_1i_2\ldots i_m} = 0$, if $i_j > mj$, for some $j$;
\item[(iii)] For each tuple $i_1,\ldots,i_m$:
\begin{equation*}\sum_{\begin{smallmatrix}\epsilon_1,\ldots,\epsilon_m\in\{0;1\}^m\\
\epsilon_1^2+\ldots + \epsilon_m^2 \ne 0
\end{smallmatrix}}p_{i_1 + m\epsilon_1,\ldots, i_m + m\epsilon_m}
= 0.\end{equation*}
\end{itemize}
\end{itemize}
\end{theor}

\begin{proof} (a) is by Theorem 4.

(b) For $m = 1$, we have $\chi(x) = (x - \lambda)$, where
$\lambda\in\Bbbk$. But each one-dimensional $\frac{d}{dt}$-invariant
subspace in $J(\lambda, 1)$ is generated by a vector function
$\overline{g}$ satisfying $\frac{d}{dt}\overline{g} =
\lambda\overline{g}$ that is by $(\alpha_1e^{\lambda t},\ldots,
\alpha_qe^{\lambda t})$, where $\alpha_i\in\Bbbk$. Hence,
$\Grass_{1}^{\frac{d}{dt}}(J(\lambda, 1))$ is a projectivization of
the linear span of such functions, i. e. it is isomorphic to
$\mathbb{P}^{q - 1}$.

(c) Before starting to prove this, we must confess that the
relations of group (ii) are in fact unnecessary, for they follow
from (i) and (iii). Although, in practice they may help to simplify
much of the group (iii) relations and to shorten their list, so we
couldn't help mentioning them. Because of this, we first prove that
the equations (ii) are satisfied in our variety, and then we show
that it is in fact given by (i) and (iii).

It is necessary to fix some notation. For $i = 1,\ldots,q$, $j =
0,\ldots, m - 1$, set $e_{ij} = (0, \ldots, 0,
\frac{1}{j!}t^je^{\lambda t}, 0, \ldots, 0)$, where the only nonzero
component is the $i$-th one. Let now
$U\in\Grass_m^{\frac{d}{dt}}(J(\lambda, m))$. Then $U =
\lspan_{\Bbbk}\left\{\sum_{i,j}\alpha^{(k)}_{ij}e_{ij}\mid k =
1,\ldots, m \right\}$, where base elements will be chosen in the
following way. Decompose $U$ into a direct sum of subspaces $U_t$
satisfying the property that each characteristic polynomial
$\chi_t(x) = x^n + c_{n-1}x^{n-1} + \ldots + c_1x + c_0$ of the
restriction of $\frac{d}{dt}$ to $U_t$ is minimal. But in this case
in a certain base the matrix of $\frac{d}{dt}|_{U_t}$ will be
written as
$${\small
\begin{pmatrix}
-c_{n - 1} & 1 & \ & \ & \ \\
-c_{n - 2} & 0 & 1 & \ & \ \\
\vdots & \ & \ddots & \ddots &\ \\
-c_{1} & \ & \ & 0 & 1 \\
-c_0 & \ & \ & \ & 0 \\
\end{pmatrix},}$$
which means that the first basic vector is a derivative of the
second, the second one is the derivative of the third and so on
until the last one which satisfies the differential equation
$\chi_t(\frac{d}{dt})g = 0$. Collecting the bases of all $U_t$, we
obtain a convenient spanning set.

In $\mathbb{P}(\bigwedge^mJ(\lambda, m))$ to the subspace $U$
associated is a line spanned by $\omega_U =
(\sum_{i,j}\alpha^{(1)}_{ij}e_{ij})\wedge\ldots\wedge
(\sum_{i,j}\alpha^{(m)}_{ij}e_{ij})$. It is easy to see that
coefficients in the decomposition $\omega_U = \sum_{l_1m + \nu_1 <
\ldots < l_mm + \nu_m}p_{l_1m + \nu_1,\ldots,l_mm +
\nu_m}e_{\nu_1l_1}\wedge\ldots\wedge e_{\nu_ml_m}$ satisfy the
relations (ii).

Let now $\omega_U =
(\sum_{i,j}\alpha^{(1)}_{ij}e_{ij})\wedge\ldots\wedge
(\sum_{i,j}\alpha^{(m)}_{ij}e_{ij}) = \sum_{L = (l_1m + \nu_1 <
\ldots < l_mm + \nu_m)}p_{L}e_{L}$ be a tensor corresponding to a
subspace $U$. One can easily check that $U$ is
$\frac{d}{dt}$-invariant if and only if $\omega_{U}$ is a relative
$\frac{d}{dt}$-invariant. On the other hand,
$\frac{d}{dt}\cdot\omega_U =
(\sum_{i,j}\alpha^{(1)}_{ij}\frac{d}{dt}e_{ij})\wedge\ldots\wedge
(\sum_{i,j}\alpha^{(m)}_{ij}\frac{d}{dt}e_{ij}) =
(\sum_{i}(\lambda\alpha^{(1)}_{i0}e_{i0} + \alpha^{(1)}_{i1}(\lambda
e_{i1} + e_{i0}) + \ldots + \alpha^{(1)}_{im}(\lambda e_{im} +
e_{i,m-1})))\wedge\ldots\wedge
(\sum_{i}(\lambda\alpha^{(m)}_{i0}e_{i0} + \alpha^{(m)}_{i1}(\lambda
e_{i1} + e_{i0}) + \ldots + \alpha^{(m)}_{im}(\lambda e_{im} +
e_{i,m-1}))) = \lambda\omega_U + \yen_U$. The last term denoted by
$\yen_U$ is to be investigated. Its summands are obtained when
$e_{i,j-1}$ are taken instead of $e_{ij}$ in the above wedge
product. Therefore, the coefficient of
$e_{\nu_1l_1}\wedge\ldots\wedge e_{\nu_ml_m}$ in this term equals
$$\sum_{\begin{smallmatrix}\epsilon_1,\ldots,\epsilon_m\in\{0;1\}^m\\
\epsilon_1^2+\ldots + \epsilon_m^2 \ne 0
\end{smallmatrix}}p_{ml_1 + \nu_1 + m\epsilon_1,\ldots, ml_m + \nu_m + m\epsilon_m}.$$

But let $l^*_1,\ldots,l^*_m$ be such a tuple that $p_{mr_1 +
\eta_1,\ldots,mr_m + \eta_m} = 0$ or not defined for all vectors
$(mr_1 + \eta_1,\ldots,mr_m + \eta_m)$ with $r_i
> l^*_i\,\forall i$, and, moreover, $p_{mr_1 + \mu_1,\ldots,mr_m + \mu_m} \ne 0$
for some $\mu_1,\ldots,\mu_m$. Then the coefficient
$e_{\mu_1r_1}\wedge e_{\mu_mr_m}$ of $\frac{d}{dt}\cdot\omega_U$
equals $\lambda^mp_{mr_1 + \mu_1,\ldots,mr_m + \mu_m}$, and thus if
$\omega_U$ is a relative $\frac{d}{dt}$-invariant, we have
$\frac{d}{dt}\cdot\omega_U = \lambda^m\omega_U$. Consequently,
 $U$ is
$\frac{d}{dt}$-invariant if and only if $\yen_U$ is zero. But we
have already shown that, rewritten in terms of Plucker coordinates
of $\omega_{\lambda}$, this condition becomes (iii). This completes
the proof of the theorem.
\end{proof}

{\bf Example 5}. For $m = q = 2$ these relations are very simple. It
is a straightforward computation to check that in this case, for
$\mu\ne\lambda$
$$\pi_{s}^{-1}((x - \lambda)(x - \mu))\cong\mathbb{P}^1\times\mathbb{P}^{1},$$ which is a non-degenerate
quadric. When the eigenvalues coincide, the fiber is
$$\pi_{s}^{-1}((x - \lambda)^2)\cong\left\{\omega = \sum_{k <
l}p_{kl}e_k\wedge e_l\in\mathbb{P}\left({\bigwedge}^2 J(\lambda,
2)\right)\left|\begin{matrix}p_{14}^2 - p_{13}p_{24} = 0,\\p_{34} = 0,
p_{14} = p_{23}\end{matrix}\right.\right\},$$ i.e. a degenerate
quadric.

Now let $Q$ be an arbitrary quiver of type $A^{(1)}_n$ (we use the
notation from the previous section for its vertices and arrows ). As
it was shown before, the layer of $\mathcal{M}^s(Q, \alpha, \zeta)$
over a point $\overline{\chi} = (\chi_1,\ldots,\chi_n)$ of the
standard categorical quotient is isomorphic to $\prod_{i =
1}^N\Grass_{\underline{r}_j}^{\Bbbk Q}(J(\lambda_j,
\underline{r}_j))$, where the notation is as in Theorem 4. On the
other hand,
$$\Grass_{\underline{r}_j}^{\Bbbk Q}(J(\lambda_j,
\underline{r}_j)) =
\left\{\left.\left(N_i\subseteq\mathcal{A}[t]\otimes\bigoplus_{j =
1}^mV_j\right)_{i = 1}^n\right|\begin{matrix}
\frac{d}{dt}(N_i)\subseteq N_i\\
a_i(N_i)\subseteq N_{i + 1}\\
(\frac{d}{dt} - \lambda_j\idd)^{\max_i{r_ij}}|_{N_i}\equiv 0
\end{matrix}\right\},$$
where $a_i: \mathcal{A}[t]\otimes\bigoplus_{j = 1}^mV_j\rightarrow
\mathcal{A}[t]\otimes\bigoplus_{j = 1}^mV_j$ acts as $\frac{d}{dt}$
on $\mathcal{A}[t]\otimes V_i$ and trivially on all other
components.

Convenient is to fix the following basis in all $J(\lambda,
\underline{r}_j)_i$: $e^{(i, \lambda)}_{qrs} =
(0,\ldots,t^se^{\lambda t},\ldots,0)\in\mathcal{A}[t]\otimes V_q$,
where the only nonzero component is the $r$-th one. Furthermore, in
$\Grass_{r_{ij}}(J(\lambda_j, \underline{r}_j)_i)$ we will be
considering Plucker's coordinates corresponding to this basis; they
will be denoted by $p^{(j)}_{k_1,\ldots,k_{\alpha_j}}$, where $k_h$
are in fact triples of indices $(q_h,r_h,s_h)$. Recall that by
Theorem 4 it is sufficient to describe the quotient for the case,
when all $\lambda_i$ coincide.

\begin{theor} Let $Q$ be a quiver of type $A^{(1)}_n$.
The Grassmannian $\Grass_{\alpha}^{\Bbbk Q}(J(\lambda, \alpha))$ is
isomorphic to the subvariety of $\prod_{i =
1}^n\mathbb{P}\left(\bigwedge^{\alpha_i}\left(J(\lambda,
\alpha)\right)_i\right)$ given by:
\begin{itemize}
\item[(a)] The equations (i) -- (iii) from Theorem 5 for each
component $P\left(\bigwedge^{\alpha_i}\left(J(\lambda,
\alpha)\right)_i\right)$;
\item[(b)] $p^{(i)}_{j_1\ldots j_{\alpha_i}}p^{(i +
1)}_{k_1\ldots k_{\alpha_{i}}} = p^{(i)}_{k_1\ldots
k_{\alpha_i}}p^{(i + 1)}_{k_1\ldots k_{\alpha_{i}}}$, for all
$(j_1,\ldots, j_{\alpha_i})$ and $(k_1,\ldots, k_{\alpha_1})$, for
$\lambda\ne 0$;
\item[(c)] $\sum_{l = 0}^{\alpha_{i + 1}}(-1)^lp^{(i)}_{j'_1\ldots j'_{\alpha_i - 1}k'_l}
p^{(i + 1)}_{k_0\ldots \widehat{k_l}\ldots k_{\alpha_{i + 1}}} = 0$,
for all tuples $j_1,\ldots,j_{\alpha_i - 1}$ (the last component in
these triples has to be less than $\alpha_i$) and $k_1,\ldots,
k_{\alpha_{i + 1}}$, where $k'_l = (q,r,s + 1)$, if $k_l = (q,r,s)$
with $q = i$, and $k'_l = k_l$ otherwise, for $\lambda = 0$.
\end{itemize}
\end{theor}

\begin{proof} 


At first, all $M_i$ are $\frac{d}{dt}$-invariant, so (i) -- (iii)
of Theorem 5 hold. In addition, there is a condition
$a_i(M_i)\subseteq M_{i + 1}$. If $\alpha_i = \alpha_{i + 1}$, this
is equivalent to the corresponding skew symmetric tensors being
proportional, i. e. to (b).

Otherwise, if $\lambda = 0$, we use a technique that may in fact be
applied in any case. Namely, let $\omega_{M_i}\in
P\left(\bigwedge^{\alpha_i}(J(\lambda, \alpha))_i\right)$ and
$\omega_{M_{i + 1}}\in P\left(\bigwedge^{\alpha_{i + 1}}(J(\lambda,
\alpha))_{i + 1}\right)$ be the tensors corresponding to $M_i$ and
$M_j$; let also $\omega'_{M_i} = a_i\cdot\omega_{M_i}$. Consider the
basis $\xi^{(i + 1)}_{qrs}$ in $\bigwedge^{\alpha_i}\left(J(\lambda,
\alpha)\right)_i$, that is adjoined to $e^{(i)}_{q,r,s}$. Then
$a_i(M_i)\subseteq M_{i + 1}$ if and only if for each tuple
$j_1,\ldots,j_{\alpha_{i} - 1}$
we have $\omega'_{M_i}(\xi^{(i + 1)}_{j_1},\ldots,\xi^{(i +
1)}_{j_{\alpha_i - 1}}, \cdot)\wedge\omega_{M_{i + 1}} = 0$. But, on
the other hand, $\langle a_i(e^{(i)}_{k}), \xi^{(i +
1)}_{j}\rangle\ne 0$ if and only if $k = j'$ (the meaning of prime
was defined in the statement of theorem), so $\omega'_{M_i}(\xi^{(i
+ 1)}_{j_1},\ldots,\xi^{(i + 1)}_{j_{\alpha_i - 1}}, \cdot) =
\sum_{l}p^{(i)}_{j'_1\ldots j'_{\alpha_{i} - 1} l'}e^{(i + 1)}_{l}$.
Consequently, $\omega'_{M_i}(\xi^{(i + 1)}_{j_1},\ldots,\xi^{(i +
1)}_{j_{\alpha_i - 1}}, \cdot)\wedge\omega_{M_{i + 1}} =
(\sum_{l}p^{(i)}_{j'_1\ldots j'_{\alpha_{i} - 1} l'}e^{(i +
1)}_{l})\wedge(\sum p^{(i + 1)}_{h_1\ldots h_{\alpha_{i + 1}}}e^{(i
+ 1)}_{h_1}\wedge\ldots\wedge e^{(i + 1)}_{h_{\alpha_{i + 1}}})$.
Having opened the brackets, we obtain that the coefficient of $e^{(i
+ 1)}_{k_0}\wedge e^{(i + 1)}_{k_1}\ldots e^{(i + 1)}_{k_{\alpha_{i
+ 1}}}$ equals $\pm\sum_{l = 0}^{\alpha_{i +
1}}(-1)^lp^{(i)}_{j'_1\ldots j'_{\alpha_i - 1}k'_l} p^{(i +
1)}_{k_0\ldots \widehat{k_l}\ldots k_{\alpha_{i + 1}}}$. But as was
mentioned before, the inclusion $a_i(M_i)\subseteq M_{i + 1}$ is
equivalent to the fact that all $\omega'_{M_i}(\xi^{(i +
1)}_{j_1},\ldots,\xi^{(i + 1)}_{j_{\alpha_i - 1}},
\cdot)\wedge\omega_{M_{i + 1}}$ are zero, which means that their
coefficients are zero. Thus we come to (c).

Theorem 6 is proved.
\end{proof}

\section{Quivers with successive cycles}


A quiver $Q$ will be called a {\it quiver with successive cycles},
if whenever two oriented cycles in $Q$ have a common vertex, they
are both powers of a certain cycle. It is easy to see that all such
quivers may be constructed through the following procedure (that
justifies our choice of terminology). Take a quiver without oriented
cycles and replace some of its vertices by oriented cycles so that
the arrows that used to start from the replaced vertex may now start
from any chosen vertex of the pasted cycle. Here is an example of a
quiver obtained through such a transformation:

\begin{center}
\includegraphics[scale=0.6]{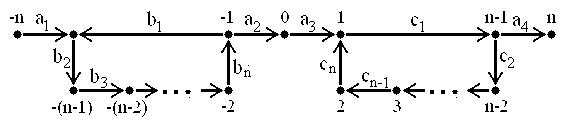}
\end{center}

Using our description of the quotient for the case of one oriented
cycle, we can generalize the technique we possess to such quivers.
As before, we denote by $\pi_s$ the natural projection
$\mathcal{M}^s(Q, \alpha, \zeta)\rightarrow\mathcal{M}(Q, \alpha,
\zeta)$.

\begin{theor} Let $Q$ be a quiver with successive cycles. Let also
$\alpha$ and $\zeta$ be two dimension vectors and $y$ be a point in
$\mathcal{M}(Q, \alpha, \zeta)$. There exists a quiver
$Q^{\spadesuit}$, a dimension vector
$\widetilde{\alpha}\in(\mathbb{Z}_{\geqslant
0})^{Q^{\spadesuit}_{0}}$, and a finite dimensional representation
$W^{\spadesuit}$ of $Q^{\spadesuit}$ such that
$\pi_s^{-1}(y)\cong\Grass^{\Bbbk
Q^{\spadesuit}}_{\widetilde{\alpha}}(W^{\spadesuit})$.
\end{theor}

\begin{proof} We begin the proof
by giving a construction of a module $J^{!!}$ such that points in
$\Grass^{\Bbbk Q}_{\alpha}(N)$ are in one-to-one correspondence with
points in the quotient space $\mathcal{M}^s(Q, \alpha, \zeta)$.

Let $\widehat{Q}$ be the quiver without oriented cycles, from which
$Q$ may be obtained through the above procedure, and $\Xi_i$ be the
set of all paths in $\widehat{Q}$ starting at $i$. For the vertices
of the oriented cycle that we place instead of the $i$-th vertex of
$\widehat{Q}$ we set the following notation:
$$\xymatrix{
i^{(0)}\ar[r]^{c_{i, 1}} & i^{(1)}\ar[r]^{c_{i, 2}}&
i^{(2)}\ar[r]^{c_{i, 3}}&\ldots\ar[r]^{c_{i, n_i - 1}} & i^{(n_i - 1)}\ar[dll]\\
&& i^{(n_i)}\ar[ull]^{c_{i, n_i}} &&}.$$ Let also $\tau_{i,j}$ be
the cycle of minimal nonzero length starting at $i^{(j)}$, if there
is any, or $e_i = e_{i^{(0)}}$, otherwise.

 As before, we consider the $\Bbbk
Q$-module $J$ with $J_i = \prod_{\sigma:i\rightsquigarrow
j}V_j^{\sigma}$. Every path $\sigma$ in $Q$ is of the form
$\sigma = B_l\tau_{i_l, j_l}^{k_l}a_{l}\ldots
B_2\tau_{i_2,j_2}^{k_2}a_2B_1\tau_{i_1, j_1}^{k_1}a_1B_0\tau_{i_0,
j_0}^{k_0}$, where $l = l(\sigma)$ is the length of $\sigma$,
$a_{j}$ are arrrows of $\widehat{Q}$, and $B_t$ are segments of
$\tau_{i_t,j_t}$, that are not oriented cycles. Let $a_0 = i$ and
$u_t = ha_t$. We then have that, for $r = 0,\ldots,n_{u_l}$,
$$J_{i^{(r)}}\cong\prod_{\Xi_i\ni\rho =
a_l\ldots a_2a_1}\bigoplus_{t = 0}^{\max(n_{u_l} - 1,
0)}K_l[x_{\rho,l}]\otimes\ldots\otimes K_0[x_{\rho, 0}]\otimes
V_{u_l^{(t)}},$$ where
$$K_q[x] := \begin{cases}
\Bbbk[[x]],\mbox{ if $l(\tau_{i_q, j_q}) > 0$},\\
\Bbbk,\mbox{ otherwise}.
\end{cases}$$
Thus we may set
$$J^{!!}_{i^{(r)}} := \bigoplus_{\Xi_i\ni\rho =
a_l\ldots a_2a_1}\bigoplus_{t = 0}^{\max(n_{u_l} - 1,
0)}\mathcal{K}_l[x_{\rho,l}]\otimes\ldots\otimes
\mathcal{K}_0[x_{\rho, 0}]\otimes V_{u_l^{(t)}},$$ where
$$\mathcal{K}_q[x] := \begin{cases}
\mathcal{A}[x],\mbox{ if $l(\tau_{i_q, j_q}) > 0$},\\
\Bbbk,\mbox{ otherwise}.
\end{cases}$$
Let $c_{1} = c_{u_l, 1}$ be the arrow belonging to $\tau_{i_l,j_l}$
that starts at $u_l$. For a pair $(M, f)\in\repq(Q, \alpha,\zeta)$
the map $\Phi^!_{(M, f)}: M\rightarrow J^{!!}$ acts as
$$\varphi_{i^{(r)}} = \bigoplus_{\Xi_i\ni\rho =
a_l\ldots a_2a_1}\!\!\!\!\!\bigoplus_{t = 0}^{\max(n_{u_l} - 1,
0)}\!\! f_{i_l^{t + 1}}c_{t}\ldots c_{1}\mathfrak{exp}(\tau_{i_l,
j_l}x_{\rho,l})\ldots
a_1B_0\mathfrak{exp}(\tau_{i_0, r}x_{\rho,0}), \eqno{(9)}$$ where
$$\mathfrak{exp}(\tau_{i, j}y) := \begin{cases}
\exp(\tau_{i, j}y),\mbox{ if $l(\tau_{i, j}) > 0$},\\
\idd,\mbox{ otherwise}.
\end{cases}$$

Combining the proofs of Lemma 2 and Lemma 6 we see that thus defined
$\varphi^!_{(M, f)}$ enjoys its usual properties, i.e. that the
statement of Lemma 6 holds in this situation.


\begin{utver} The map $\Phi^!: \repq^s(Q, \alpha, \zeta)\rightarrow\IHom_{\alpha}^{\Bbbk
Q}(J^{!!})$, $(M, f)\mapsto\Phi^!_{(M,f)}$ is a bijection.
\end{utver}

\begin{proof} Let $U\subseteq J^{!!}$ be a subrepresentation with
dimension vector $\alpha$. We need a pair $(M, f)$ such that
$\TrueIm(\Phi^!_{(M,f)}) = U$. Take $M = U$. The maps $f_{i^{(j)}}$
are then reconstructed as compositions of the projections
$U_{i^{(j)}}\rightarrow \mathcal{K}_0[x_{e_i, 0}]\otimes
V_{i^{(j)}}$ with evaluation at 0. It is now only left to show that
$\Phi^!_{(U,f)}$ is the natural inclusion $U\hookrightarrow J^{!!}$.


The proof is by ``downward induction over $\widehat{Q}$''. If
$i\in\widehat{Q}$ is a sink, we apply Proposition 5. For an
arbitrary vertex $i$, fix a basis $G_1,\ldots,G_{d_j}$ in each
$U_{i^{(j)}}$ and observe that, as in (5),
$$\left(G_1(x_{e_i,0}),\ldots G_{d_j}(x_{e_i,0})\right) =
\left(G_1(0),\ldots G_{d_j}(0)\right)\exp(\tau_{i,j}x_{e_i,0}),$$
where $G_{p}(x_{e_i},0)$ are in fact functions in $x_{\rho,t}$, for $\rho\ne e_i$. We need to prove that
the right hand side of this equality coincides with (9).
Recall that $\bigoplus_{a: i\rightarrow p} aB_{0,a}$, where $B_{0,a}$ is the shortest segment of $\tau_{i,j}$
linking $i^{j}$ with $ta$, acts as evaluation at
$x_{e_i, 0} = 0$. So, $\left(G_1(0),\ldots G_{d_j}(0)\right)$ consists of horizontal blocks
representing the bases of $aB_{0,a}(U_{i^{j}})$, where $a$ are arrows starting at
vertices of $\tau_{i,j}$. By the induction hypothesis all these blocks
are of the form (9). Thus the claim follows. 
\end{proof}

\begin{corol} Points of the quotient space $\mathcal{M}^s(Q, \alpha,
\zeta)$ are in one-to-one correspondence with points of the
Grassmannian of submodules $\Grass_{\alpha}^{\Bbbk Q}(J^{!!})$.
\end{corol}

This, however, rather characterizes $\Grass_{\alpha}^{\Bbbk
Q}(J^{!!})$, than the quotient space. So, as we did in Section 4, we
restrict our attention to the fibers of $\pi_s$.

Recall that the algebra $\Bbbk[\repq(Q, \alpha, \zeta)]$ is
generated by coefficients of characteristic polynomials of all
oriented cycles in $Q$. So we can consider $y$ as a tuple
$\left\{\widehat{\chi}_{i,j}\mid i\in\widehat{Q}_0,\, j = 0,\ldots,
n_i\right\}$, where
$$\widehat{\chi}_{i,j} := \begin{cases}
\mbox{the characteristic polynomial of $\tau_{i, j}$, if $l(\tau(i, j) > 0$},\\
0,\mbox{ otherwise}.
\end{cases}$$

We are now ready to prove the theorem. Introduce the module $W$ with
$$W_{i^{(r)}} := \bigoplus_{\Xi_i\ni\rho =
a_l\ldots a_2a_1}\left(\bigotimes_{t =
0}^l\mathcal{K}_t[x_{\rho,t}]\otimes\left(\bigoplus_{t =
0}^{\max(n_{u_l} - 1, 0)}
V_{u_l^{(t)}}\right)\right)_{\begin{smallmatrix}
\widehat{\chi}_{i_t,j_t}\left(\frac{\partial}{\partial x_{\rho, t}}\right),\\
t = 0,\ldots,l
\end{smallmatrix}},$$
where $j_0 = r$. It is clearly finite dimensional. We claim that for
a stable pair $(M, f)\in\pi_s^{-1}(y)$ the image of
$\varphi^{!}_{(M, f)}$ lies in $W$.

We prove this by ``downward induction over $\widehat{Q}$''. Let $i$
be a sink in $\widehat{Q}$. If this vertex is replaced by an
oriented cycle, we may use the results of Section 4, otherwise the
claim is trivial.

Now, consider an arbitrary vertex $i$. If it does not belong to an
oriented cycle, recall that by Proposition 3 we have, for every
submodule $U$ of $J^{!!}$, $U_{i^{(0)}}\subseteq
V_{i^{(0)}}\oplus\bigoplus_{a:i^{(0)}\rightarrow
j^{(r)}}U_{j^{(r)}}\subseteq
V_{i^{(0)}}\oplus\bigoplus_{a:i^{(0)}\rightarrow
j^{(r)}}W_{j^{(r)}}$, but by definition of $\widehat{\chi}_{i, 0}$
this clearly lies in $W_{i^{(0)}}$. Otherwise, we use (9) with
$\mathfrak{exp}$ the usual matrix exponents to understand that all
$\widehat{\chi}_{i_t,j_t}(\frac{\partial}{\partial x_{\rho, t}})$
annulate $W_{i^{(r)}}$, for $r = 0,\ldots,n_i$.

Although $\Phi^!: (M, f)\mapsto\Phi^{(M,f)}$ induces an inclusion
$\pi^{-1}_{GL(\alpha)}(y)\cap\repq^s(Q,\alpha,\zeta)\hookrightarrow\IHom_{\alpha}^{\Bbbk
Q}(W)$, in general it is not surjective. For example, let $Q$ be a
Jordan quiver consisting of a single loop $a$. Take $\alpha = \zeta
= 2$ and $\chi_a(x) = (x - 1)(x - 2)$. Then the corresponding $W$
contains a subrepresentation $U$ spanned by
$\left(\begin{smallmatrix}e^{t}\\ 0\end{smallmatrix}\right)$ and
$\left(\begin{smallmatrix}0\\e^{t}\end{smallmatrix}\right)$, which
is not in $\pi^{-1}_{GL(2)}(\chi_a)$. Hence we need to refine our
construction.

Let $\lambda^{(i)}_1,\ldots,\lambda^{(i)}_{N_{i}}$ be all different
roots of $\widehat{\chi}_{i,j}$, for $j = 1,\ldots,n_i$. Let also
$r^{(i,j)}_{l}$ be the multiplicity of $\lambda^{(i)}_l$ as a root
of $\widehat{\chi}_{i,j}$ and $\underline{r}^{(i)}_l =
(r^{(i,0)}_l,\ldots,r^{(i,n_i)}_l)$. If $\widehat{\chi}_{i,j}$ is a
zero polynomial, we set $N_{i} = 1$, $\lambda^{(i)}_1 = 0$, and
$r^{(i,0)}_1 = \alpha_{i^{(j)}}$. For $j = 0,\ldots,n_i$, define
$J(\lambda^{(i)}_l, \underline{r}^{(i)}_l)_{i^{(j)}}$ as the
subspace of $W_{i^{(j)}}$ generated by all solutions of the
differential equation $(\frac{\partial}{\partial x_{e_i, 0}} -
\lambda^{(i, j)}_l\idd)^{\max_qr^{(i,q)}_{l}}\overline{g} = 0$.
Note, that usually these subspaces do not form a subrepresentation,
like they used to in Section 4. Since $W_{i^{(j)}} = \bigoplus_{l =
1}^{N_i}\left(W_{i^{(j)}}\cap J(\lambda^{(i)}_l,
\underline{r}^{(i)}_l)_{i^{(j)}}\right)$, there are natural
projections $\mathbf{p}^{(i,j)}_l: W_{i^{(j)}}\twoheadrightarrow
W_{i^{(j)}}\cap J(\lambda^{(i)}_l, \underline{r}^{(i)}_l)_{i^{(j)}}$
along $\bigoplus_{l' \ne l}^{N_i}\left(W_{i^{(j)}}\cap
J(\lambda^{(i)}_{l'}, \underline{r}^{(i)}_{l'})_{i^{(j)}}\right)$.
The same arguments as we used to prove Lemma 9, show that a stable
pair $(M,f)$ lies in $\pi^{-1}_{GL(\alpha)}(y)$ if and only if the
dimensions of $\mathbf{p}^{(i,j)}_l(\Phi^!_{(M,f)}(M)_{i^{(j)}})$
equal $r^{(i,j)}_{l}$, for all $i$, $j$ and $l$. Hence, $\Phi^!$,
induces an inclusion of $\pi^{-1}_{GL(\alpha)}(y)\cap\repq^s(Q,
\alpha, \zeta)$ into the subvariety $\widetilde{\mathrm{IH}}(W)$ in
$$\prod_{\begin{smallmatrix}i^{(j)}\in\Bbbk Q_0,\\
 l = 1,\ldots, N_i
\end{smallmatrix}}\IHom_{r^{(i,j)}_l}\left(W_{i^{(j)}}\cap J(\lambda^{(i)}_l, \underline{r}^{(i)}_l)_{i^{(j)}}\right),$$
consisting of those tuples of injections, whose images give a
$Q$-subrepresentation of $W$. Using Lemma 4 one may deduce that this
inclusion is in fact a $GL(\alpha)$-invariant isomorphism of
algebraic varieties. Therefore,
$\pi_s^{-1}(y)\cong\widetilde{\mathrm{IH}}(W)/\!\!/GL(\alpha)$.
However, the latter does not yet look like a Grassmannian of
submodules.

We now construct a quiver $Q^{\spadesuit}$, setting
$$Q^{\spadesuit}_0 = \left\{i^{(j,l)}\mid i^{(j)}\in Q_0, l =
1,\ldots, N_i\right\}$$ and $$Q^{\spadesuit}_1 = \left\{c_{i,j,l}:
i^{(j,l)}\rightarrow i^{(j+1, l)}\mid j = 0,\ldots, n_i; l =
1,\ldots, N\right\}\cup\\$$ $$\cup\left\{a^{(l_1, l_2)}: i_1^{(j_1,
l_1)}\rightarrow i_2^{(j_2, l_2)}\mid a: i_1^{j_1}\rightarrow
i_2^{j_2}\right\}.$$ Observe that $W$ may be considered as a
representation $W^{\spadesuit}$ of $Q^{\spadesuit}$ with
$W^{\spadesuit}_{i^{(j,l)}} = W_{i^{(j)}}\cap J(\lambda^{(i)}_l,
\underline{r}^{(i)}_l)_{i^{(j)}}$. The dimension vector of
$W^{\spadesuit}$ is $\widetilde{\alpha}$ with
$\widetilde{\alpha}_{i^{(j,l)}} = r^{(i,j)}_l$. It is now
straightforward to check that
$\widetilde{\mathrm{IH}}(W)/\!\!/GL(\alpha)$ is isomorphic to the
Grassmannian $\Grass_{\widetilde{\alpha}}^{\Bbbk
Q^{\spadesuit}}(W^{\spadesuit})$. This finishes the proof.
\end{proof}

{\bf Remark.} For $Q = A^{(1)}_{n-1}$ the quiver $Q^{\spadesuit}$ is
non-connected. In fact, it is a union of $N$ copies of $Q$, where
$N$ is as in Theorem 4. So, a $\Bbbk Q^{\spadesuit}$-module is a
$N$-tuple of representations of $Q$. In particular, $W^{\spadesuit}$
splits as $\bigoplus_jJ(\lambda_j, \underline{r}_j)$. So,
$\Grass_{\widetilde{\alpha}}^{\Bbbk
Q^{\spadesuit}}(W^{\spadesuit})\cong\prod_j\Grass_{\underline{r}_j}^{\Bbbk
Q}(J(\lambda_j, \underline{r}_j))$.

\medskip

To illustrate both the theorem and the proof, we give an example.

{\bf Example 6}. Let $\widehat{Q} = \xymatrix{
1\ar@<1ex>[r]^{a'_1}\ar@<-1ex>[r]_{a^{'{'}}_1} & 2\ar[r]^{a_2} &
3\ar[r]^{a_3} & 4\ar[r]^{a_4} & 5}.$ Replace its second and fourth
vertices by loops $b$ and $c$ respectively and denote the
quiver obtained by $Q$. 
Fix two dimension vectors $\alpha$ and $\zeta$ and construct our
usual $J^{!!}$. Then $J^{!!}_5 = V_5$, $J^{!!}_{4} =
\mathcal{A}[x]\otimes V_4\oplus \mathcal{A}[x]\otimes V_5$,
$J^{!!}_{3} = V_3\oplus\mathcal{A}[x]\otimes V_4\oplus
\mathcal{A}[x]\otimes V_5$, $J^{!!}_{2} = \mathcal{A}[y]\otimes
V_2\oplus \mathcal{A}[y]\otimes V_3\oplus\mathcal{A}[x, y]\otimes
V_4\oplus \mathcal{A}[x, y]\otimes V_5$, and $J^{!!}_{1} =
(V_1\oplus \mathcal{A}[y_1]\otimes V_2\oplus \mathcal{A}[y_1]\otimes
V_3\oplus\mathcal{A}[x, y_1]\otimes V_4\oplus \mathcal{A}[x,
y_1]\otimes V_5) \oplus (V_1\oplus \mathcal{A}[y_2]\otimes V_2\oplus
\mathcal{A}[y_2]\otimes V_3\oplus\mathcal{A}[x, y_2]\otimes
V_4\oplus \mathcal{A}[x, y_2]\otimes V_5)$ with maps between them as
follows: $a'_1 =
\Box_y\circ((0\oplus\idd\oplus\idd\oplus\idd\oplus\idd)\oplus(0\oplus0\oplus0\oplus0\oplus0))$,
$a'{'}_1 =
\Box_y\circ((0\oplus0\oplus0\oplus0\oplus0)\oplus(0\oplus\idd\oplus\idd\oplus\idd\oplus\idd))$,
$b = \frac{\partial}{\partial y}$, $a_2 =
0\oplus\gamma_0\oplus\gamma_0\oplus\gamma_0$, $a_3 =
0\oplus\idd\oplus\idd$, $c = \frac{d}{dx}$, and $a_4 =
0\oplus\delta_{0}$, where $\delta_{0}: G(x) \mapsto G(0)$,
$\gamma_0: F(x, y)\mapsto F(x, 0)$ and $\Box_y(F(x, y_i)) := F(x,
y)$.

Our purpose is now to show that fibers of
$\mathcal{M}^s(Q,\alpha,\zeta)$ over points of categorical quotient
are realized as a $\Bbbk Q^{\spadesuit}$-module Grassmannians
$\Grass_{\widetilde{\alpha}}^{\Bbbk Q^{\spadesuit}}(W^{\spadesuit})$
for some quiver $Q^{\spadesuit}$, dimension vector
$\widetilde{\alpha}$, and finite dimensional module
$W^{\spadesuit}$. The categorical quotient
$\repq(Q,\alpha)/\!\!/GL(\alpha)$ is isomorphic to
$\mathbb{A}^{\alpha_2}\times\mathbb{A}^{\alpha_4}$, and its points
may be viewed as pairs $\overline{\chi} = (\chi_b, \chi_c)$ of
characteristic polynomials of $b$ and $c$. Having fixed such a pair,
we construct the required modules $W$ by ``downward induction''.

At first, we need a finite dimensional module $W$ such that
$\Phi^!_{(M, f)}(M)$ lies in $W$, for each
$M\in\pi_{GL(\alpha)}^{-1}(\overline{\chi})$. As usually, let
$\lambda^{(b)}_1,\ldots,\lambda^{(b)}_{N_b}$ and
$\lambda^{(c)}_1,\ldots,\lambda^{(c)}_{N_c}$ be all different roots
of $\chi_b$ and $\chi_c$ respectively, with $r^{(b)}_{i}$ and
$r^{(c)}_i$ their multiplicities. 

We set $W_5 = J_5$, as it is already finite dimensional. Further,
since $U = (U_i)_{i = 1}^4$ is in $\Grass^A_{\alpha}(J^{!!})$, we
have that $U_3\subseteq 
\bigoplus_{j = 1}^{N_c}J(\lambda^{(c)}_j, r^{(c)}_j) = : W_4$. We
can also observe that $U_2\subseteq\bigoplus_{j =
1}^{N_b}J(\lambda^{(b)}_j, r^{(b)}_j)
$, although it tells nothing about $U_3$ or $U_1$ and, moreover,
$J(\lambda^{(b)}_j, r^{(b)}_j)$ are not finite dimensional, so we
should anyway continue our investigation. But using Proposition 3,
we can state that $U_3\subseteq V_3\oplus U_4$, and therefore
$U_3\subseteq V_3\oplus\bigoplus_{j = 1}^{N_b}J(\lambda^{(b)}_j,
r^{(b)}_j) =: W_3$. To determine $W_2$ we should recall that $U$ is
the image of $\Phi^{!}(M,f)$ for a pair $(M, f)\in\repq^s(Q, \alpha,
\zeta)$, so that
$$\forall\,u\in U_2\,\exists\,v\in M:u = (\Phi^{!}(M, f))_2(v) = $$
$$ =
(f_2\oplus f_3a_2\oplus f_4\exp(xc)a_3a_2 \oplus
f_5a_4\exp(xc)a_3a_2)\exp(yb)v.$$ Thus, $a_3a_2u =
(f_4\exp(xc)a_3a_2 \oplus f_5a_4\exp(xc)a_3a_2)v\in U_4\subseteq
(J^{!!}_3)_{\chi_c(\frac{d}{dx})}$. This shows that $U_2\subseteq
\bigoplus_{j, l}(J(\lambda^{(b)}_{j}, r^{(b)}_{j})\cap
J(\lambda^{(c)}_{l}, r^{(c)}_{l})) =: W_2$. Finally, $U_1\subseteq
V_1\oplus \left(\dag_{a_1}(W_2) + \dag_{a_2}(W_2)\right) =: W_1$.
So, we have found a finite dimensional submodule $W\subseteq J^{!!}$
such that $\Phi^!_{(M, f)}(M)\subseteq W$, for each
$M\in\pi_{GL(\alpha)}^{-1}(\overline{\chi})$. In a sense $W$ is the
solution space of the system $\chi_b(\frac{\partial}{\partial x})F =
\chi_c(\frac{\partial}{\partial y})F = 0$, though this notation is
rather abusive.

Now, for $j = 1,\ldots, N_b$ and $l = 1,\ldots,N_c$ consider the
natural projections $p^{(b)}_{j}: W_2\twoheadrightarrow W_2\cap
J(\lambda^{(b)}, r^{(b)}_j)_2$ and $p^{(c)}_{l}:
W_4\twoheadrightarrow W_4\cap J(\lambda^{(c)}, r^{(c)}_l)_4$. Note
that a stable pair $(M, f)$ lies in $\pi_{GL(\alpha)}^{-1}(\chi_b,
\chi_c)$ if and only if $b$ and $c$ act on $M_2$ and $M_4$ with
characteristic polynomials $\chi_b$ and $\chi_c$ respectively. This
is clearly equivalent to the dimensions of $p^{(b)}_{j}(\Phi^!_{(M,
f)}(M)_2)$ and $p^{(c)}_{l}(\Phi^!_{(M, f)}(M)_4)$ being equal to
$r^{(b)}_j$ and $r^{(c)}_l$ respectively. From this we can deduce
that the map $\Phi^!: (M, f)\mapsto\Phi^!_{(M, f)}$ induces an
isomorphism of $\mathcal{M}^s(Q, \alpha, \zeta)$ with the subvariety
in $\Grass_{\alpha_1}(W_1)\times\prod_{j =
1}^{N_b}\Grass_{r^{(b)}_j}(W_2\cap J(\lambda^{(b)}_j,
r^{(b)}_j))\times\Grass_{\alpha_3}(W_3)\times\prod_{l =
1}^{N_c}\Grass_{r^{(c)}_l}(W_4\cap J(\lambda^{(c)}_l,
r^{(c)}_l))\times\Grass_{\alpha_5}(W_5)$, given by the condition
that the tuple of subspaces is a subrepresentation of $W$ (i.e., we
identify $(U_{2_j})_{j = 1}^{N_b}\in\prod_{j =
1}^{N_b}\Grass_{r^{(b)}_j}(W_2\cap J(\lambda^{(b)}_j, r^{(b)}_j))$
with the subspace $\bigoplus_jU_{2_j}$ in $W_2$).

Finally, we obtain a realization of the quotient space as a
Grassmannian of subrepresentations. Construct the quiver
$Q^{\spadesuit}$. It has vertices $Q^{\spadesuit}_0 = \{1, 2_1,
2_2,\ldots, 2_{N_b}, 3, 4_1, 4_2,\ldots, 4_{N_c}, 5\}$ and arrows
$a_{1j}', a_{1j}'{'}: 1\rightarrow 2_j$, $b_j: 2_j\rightarrow 2_j$,
$a_{2j}: 2_j\rightarrow 3$, $a_{3l}: 3\rightarrow 4_l$, $c_{l}:
4_l\rightarrow 4_l$, and $a_{4l}: 4_l\rightarrow 5$, for $j =
1,\ldots, N_b$, $l = 1,\ldots,N_c$. Thus, $Q^{\spadesuit}$ is
$$\xymatrix{& 2_1\ar@(ur, ul)[]_{b_1}\ar[dr]^{a_{21}} & & 4_1\ar[dr]^{a_{41}}\ar@(ur, ul)[]_{c_1} &\\
1\ar@/^1pc/[ur]^{a'_{11}}\ar@/_1pc/[ur]^{a'{'}_{11}}
\ar@/^1pc/[dr]_{a'_{1N_b}}\ar@/_1pc/[dr]_{a'{'}_{1N_b}} & \vdots &
3\ar[ur]^{a_{31}}
\ar[dr]_{a_{3N_c}} & \vdots & 5 \\
& 2_{N_b}\ar@(dr, dl)[]^{b_{N_b}}\ar[ur]_{a_{2N_b}} & & 4_{N_c}
\ar@(dr, dl)[]^{c_{N_c}}\ar[ur]_{a_{4N_c}} &}.$$ Since $W_2 =
\bigoplus_j(W_2\cap J(\lambda^{(b)}_j))$ and $W_4 =
\bigoplus_l(W_4\cap J(\lambda^{(c)}_l))$, $W$ may be regarded as a
$\Bbbk Q^{\spadesuit}$-module $W^{\spadesuit}$ with
$W^{\spadesuit}_{2_j} = W_2\cap J(\lambda^{(b)}_j, r^{(b)}_j)$ and
$W^{\spadesuit}_{4_l} = W_4\cap J(\lambda^{(c)}_l, r^{(c)}_l)$.
Similarly, the image of $\Phi^!_{(M, f)}$ followed by the
projections $p^{(b)}_j$ and $p^{(c)}_l$ is a representation of
$Q^{\spadesuit}$. Conversely, consider a
$Q^{\spadesuit}$-subrepresentation $\widetilde{M}$ of
$W^{\spadesuit}$ with dimension vector $\widetilde{\alpha}$, where
$\widetilde{\alpha}_i = \alpha_i$, for $i\in\{1, 3, 5\}$,
$\widetilde{\alpha}_{2_j} = r^{(b)}_j$, and
$\widetilde{\alpha}_{4_l} = r^{(c)}_l$. Taking $M_i =
\bigoplus{M_{i_j}}$, for $i = 2,4$, we obtain a $\Bbbk Q$-submodule
$M$ of $W$ with dimension vector $\alpha$ that clearly belongs to
$\pi_{GL(\alpha)}^{-1}(\chi_b, \chi_c)$. Together with Lemma 4 this
yields the desired isomorphism $\pi_s^{-1}(\chi_b,
\chi_c)\cong\Grass_{\widetilde{\alpha}}^{\Bbbk
Q^{\spadesuit}}(W^{\spadesuit})$.

\medskip

As a last remark in this section, we note that in the proof of
Proposition 2 we never used the fact that the algebra $A$ is finite
dimensional. We thus obtain

\begin{utver} Let $Q$ be a quiver, $M$ be its representation with dimension vector
$\alpha$ and $V$ be a graded vector space with dimension vector
$\zeta$. There exists a map $f: M\rightarrow V$ making a pair $(M,
f)$ stable if and only if $\mathrm{soc}\,M$ may be embedded in
$\mathrm{soc}\,J^{!!}$.
\end{utver}

The interpretation of this condition is now a bit more
sophisticated, for $S(i)$, $i\in Q$ are not the only simple $\Bbbk
Q$-modules. In fact, each oriented cycle $\tau$ of minimal length
defines a family of simple modules $S(\tau, \lambda)$, for
$\lambda\in\Bbbk\backslash\{0\}$. These are the representations with
$S(\tau, \lambda)_p = 0$, when $p$ does not belong to $\tau$, and
$S(\tau, \lambda)_p = \Bbbk$, otherwise, with $\tau$ acting by
$\lambda\idd$. So, $\mathrm{soc}\,M$ may be embedded in
$\mathrm{soc}\,J^{!!}$ if and only if the multiplicity of each
$S(i)$ and $S(\tau, \lambda)$ in $\mathrm{soc}\,M$ is not greater
than the one in $\mathrm{soc}\,J^{!!}$. And in
$\mathrm{soc}\,J^{!!}$ the multiplicity of $S(i)$ equals $\zeta_i$,
while the multiplicity of $S(\tau, \lambda)$ equals $\sum_i\zeta_i$,
where $i$ runs through the vertices of $\tau$.

\end{document}